
\input amstex
\documentstyle{amsppt}
\magnification=\magstep1
\TagsOnLeft
\NoBlackBoxes 
\pageheight{9truein}
\pagewidth{6.3truein}

\def \BarnesTwo{2}
\def \Bosma{3}
\def \CasBk{4}
\def \DavOne{5}
\def \DavTwo{6}
\def \Enn{7}
\def \fuk{8}
\def \Grace{9}
\def \Khin{10}
\def \papertwo{11}    
\def \paperone{12}    
\def \Perron{13}
\def \paperthree{14}    
\def \paperfour{15}   
\def \Pitman{16}
\def \RockBk{17}


\def\ve{{\varepsilon}}

\define\Z {\bold Z}
\define\R {\bold R}
\define\N {\bold N}
\define\Zp {\bold Z_p}
\define\Q {\bold Q}


\topmatter
\title  Lower  Bounds on the Two-sided Inhomogeneous Approximation Constant  \endtitle 

\leftheadtext{CHRISTOPHER G. PINNER}
\author{ Christopher G. Pinner  }\endauthor
\rightheadtext{Bounds on the Inhomogeneous Approximation Constant}
 
\date May 1998, Revised January 2000 \enddate
 
\abstract{ For an irrational real number $\alpha$ and 
real number $\gamma$ one defines the two-sided approximation constant
$$ M(\alpha,\gamma):=\liminf_{|n|\rightarrow \infty} |n| ||n\alpha -\gamma||.$$
We are interested here in the case of worst inhomogeneous approximation for $\alpha$
$$\rho (\alpha):=\sup_{\gamma \not\in \Z+\alpha\Z} M(\alpha,\gamma),$$
and in how $\rho (\alpha)$ is affected by the liminf of the partial quotients 
in the nearest integer continued fraction expansion of $\alpha $;
$$ \alpha=
a_0+\frac{\ve_1}{a_1 + {\displaystyle \frac{\ve_2}{a_2+{\displaystyle \frac{\ve_3}{a
_3+\cdots }}}}},\;\;\; \ve_{i}=\pm 1. $$
In particular, setting  $R:=\liminf_{i\rightarrow \infty} a_{i}$, we obtain the  optimal  lower bound
$$ \rho (\alpha ) \geq \frac{1}{4}\left(1-\frac{3}{R} + O\left(\frac{1}{R^2}\right)\right),  $$
complementing an old  optimal upper bound of Fukasawa
$$ \rho (\alpha) \leq \frac{1}{4}\left( 1-\frac{2}{R^2} +O \left(\frac{1}{R^4}\right)\right). $$
For any $\alpha$ we show how to construct a $\gamma_*$ with
$$ M(\alpha,\gamma_*)\geq \frac{1}{25.159...}. $$ 
} \endabstract

\address{ Mathematics and Computer Science, 
University of Northern British Columbia,
3333 University Way, Prince George, B.C., Canada V2N 4Z9.
{\tt pinner\@math.ksu.edu}} 
\endaddress

\endtopmatter
\document
\baselineskip=18pt 
\head 1. Introduction  \endhead

For an irrational real  number $\alpha$ and real number $\gamma$ one defines the {\it two-sided inhomogeneous approximation constant}
$$  M(\alpha,\gamma):=\liminf_{|n|\rightarrow \infty} |n| ||n\alpha -\gamma ||.\tag 1.1 $$
Fixing $\alpha$ and varying $\gamma$ gives an
{\it inhomogeneous  spectrum} for  $\alpha$ 
$$  {\bold L}(\alpha):=\{ M(\alpha,\gamma) \; \; : \;\;\gamma \in \R,\;\; \gamma \not\in \Z+\alpha \Z \}, \tag 1.2 $$
akin (awkward reciprocals aside) to  the
traditional Lagrange spectrum  
$${\bold L}:=\{ M(\alpha,0)^{-1} : \alpha \in \R\}. \tag1.3 $$
For a  $\gamma$ not of the 
form $m+n\alpha$, $n,m\in \Z$,  we showed in [\paperone]  how to
express $  M(\alpha,\gamma)$
in terms of a semi-regular continued fraction expansion
of $\alpha$
$$\align  \alpha & =[a_0; \ve_1a_1,\ve_2a_2,\ve_3a_3,....],\;\;\;\;\ve_{i}=\pm 1, \\
  & = a_0+\frac{\ve_1}{a_1 + {\displaystyle \frac{\ve_2}{a_2+{\displaystyle \frac{\ve_3}{a
_3+\cdots }}}}}, \tag1.4 \endalign     $$
(we use Bosma's notation [\Bosma]) and a  corresponding alpha-expansion  of $\{\gamma\}$
$$  \{ \gamma \}=\sum_{i=1}^{\infty} c_i D_{i-1},\;\;\; D_i:= |q_i\alpha -p_i|, \tag 1.5  $$
where $p_i/q_i$ are the  convergents arising from the truncations of the expansion of $\alpha$.
As usual $||x||$ and $\{x\}$ denote the distance from $x$ to the nearest 
integer and the fractional part of $x$ respectively. By working with $\pm \alpha$ 
(mod 1) we can of course assume that 
$$0<\alpha <\frac{1}{2},\;\;\;\;\;\; 0\leq \gamma <1.\tag 1.6 $$
We are interested here in the 
largest point of  the spectrum (1.2), 
$$\rho(\alpha):= \sup_{\gamma \not\in \Z+\alpha \Z} M(\alpha,\gamma), \tag 1.7 $$
and in particular how $\rho(\alpha)$ is affected by the size
of the partial quotients in the expansion of $\alpha$. From a well-known  result of Minkowski
$$\rho (\alpha)\leq \frac{1}{4}. \tag 1.8 $$
  Grace [\Grace] 
 showed that this bound is achievable
(see also Cassels [\CasBk, Thm IIB, p.48], Rockett-Szusz [\RockBk]
and related results of Barnes [\BarnesTwo]).  Equality though  requires
unbounded partial quotients and for badly approximable $\alpha$
Khinchin [\Khin] obtained the improvement
$$ \rho (\alpha)\leq \frac{1}{4}\sqrt{1-4M^2(\alpha,0)}, \tag 1.9 $$
with equality in (1.9) for
the golden ratio or when $\alpha$ has  period 
$$ \ve_ia_i =+R, \hbox{ $R$ even. } \tag 1.10 $$
However while $M(\alpha,0)$ is  controlled primarily  by
the size of the limsup of the partial quotients,
$\rho (\alpha)$ depends more naturally  on the liminf of the partial quotients.
It will be convenient here to use the {\it nearest integer continued fraction
expansion}, rather
than the usual {\it regular continued fraction expansion} with all $\ve_i=+1$.
That is to say  the partial quotients $a_i$ are generated iteratively by
taking the nearest integer in the continued fraction
algorithm, rather than always rounding down as in the regular 
continued fraction (see (2.3)).
When (1.4) is the nearest integer expansion and
$$ R:=\liminf_{i\rightarrow \infty}a_{i}, \tag 1.11 $$
Fukasawa [\fuk,I,II] showed that
$$\rho (\alpha)=\frac{1}{4} \hbox{ iff }  R=\infty. \tag 1.12 $$
For finite $R$ he gave a number of upper bounds. In particular 
$$ \rho (\alpha ) \leq  \frac{1}{4}\frac{R}{\sqrt{R^2+4}} =\frac
{1}{4}\left( 1-\frac{2}{R^2}+O\left(\frac{1}{R^4}\right)\right),  \hbox{ if $R$ is even,}  \tag1.13$$  
with equality when $\alpha$ has eventual period (1.10) (and an improvement otherwise), and
$$ \rho (\alpha ) \leq 
\frac{1}{2}\frac{R-1}{\sqrt{(R+1)^2-4}+\sqrt{(R-1)^2+4},}\leq \frac{1}{4}\left(1-\frac{1}{R}\right),  \hbox{ if $R\geq 5 $ is odd.}  \tag1.14 $$
Although asymptotically optimal, further  refinements  are possible 
when more is know about the signs of the $\ve_i$
(see Fukasawa [\fuk,IV]). Indeed if the nearest integer expansion coincides with the {\it  negative 
continued fraction expansion} (i.e. the $\ve_i=-1$ for all $i\geq 2$), then 
$$\rho(\alpha)\leq \frac{1}{4}\left(1-\frac{1}{R}\right), \hbox{ if $R\geq 3 $,} \tag1.15 $$ 
for even as well as 
odd $R$ (again best possible, see [\paperthree]).
Notice that (1.12) would not be true if we replaced the nearest integer by
the regular or negative 
continued fraction expansion
(from the correspondence between the regular expansion and negative expansion;
$$ [0; a_1,a_2,a_3,a_4,a_5,\ldots ]=[0; (a_1+1),\underbrace{-2,...,-2}_{a_2-1},-(a_3+2),
\underbrace{-2,...,-2}_{a_4-1},-(a_5+2),...] \tag 1.16 $$
a nearest integer expansion with $\rho (\alpha)=\frac{1}{4}$ and all $\ve_i=-1$ would
have every other $a_i=1$ in its regular expansion, and one with all
$\ve_i=1$ would have arbitrarily long strings of $a_i=2$ in its negative 
expansion).

We are interested here in obtaining corresponding lower bounds for $\rho (\alpha)$ in terms of $R$.
We showed in [\paperone]  that when all but finitely many of 
the partial quotients of $\alpha$ 
are  even 
(with additional restrictions imposed on occurrences of $\ve_ia_i=-2$) the 
value of $\rho(\alpha)$ is achieved by taking
$$ \gamma^{*} =\gamma^*(\alpha):=\sum_{i=1}^{\infty} \frac{1}{2}(a_i-1+\ve_{i})D_{i-1}. \tag 1.17 $$
If there are odd partial quotients we can no longer take the integer coefficients $c_i$ to
be all $\frac{1}{2}(a_i-1+\ve_{i})$ but it seems reasonable to take $c_i$ to be  the closest integer possible
to this. Hence if $a_{n_{i}}$ denotes the $i$th odd partial quotient we set
$$ \gamma_*=\gamma_*(\alpha):=\sum_{i=1}^{\infty} \frac{1}{2}(a_i-1+\ve_{i}+\lambda_{i})D_{i-1}, \;\;\;\lambda_{i}:=\cases 0, & \hbox{ if $a_i$ is even, } \\ (-1)^j, & 
\hbox{ if $a_i=a_{n_{j}}$. } \endcases \tag1.18 $$
The value of $M(\alpha,\gamma_*)$ then gives us a lower bound on $\rho (\alpha)$:   

\proclaim{Theorem 1} 
Suppose that (1.4) gives the nearest integer continued fraction expansion
of $\alpha$. Then with $R$ as in (1.11)
and $\gamma_*$ as in (1.18) we have
$$ \rho (\alpha) \geq  M(\alpha,\gamma_*) \geq \cases {\displaystyle \frac{1}{2(\sqrt{13}+10)}=\frac{1}{27.211...} }, & \hbox{ if $R=2$, } \\  \\  
{\displaystyle \frac{1}{20.487...} }, & \hbox{ if $R=3$, } \\  \\ {\displaystyle
    \frac{1}{4} \left( 1 -\frac{3}{R}+\frac{4}{R^2}-\frac{7}{R^3}\right)}, & \hbox{ if $R\geq 4$ is even,} \\ \\
 {\displaystyle
    \frac{1}{4} \left( 1 -\frac{3}{R}+\frac{3}{R^2}-\frac{6}{R^3}\right)}, & \hbox{ if $R\geq 5$ is odd.} \endcases \tag 1.19 $$
\endproclaim
 
Theorem 1 will follow from the more precise Theorem 4 bound 
given in  Section 3.
Our result is asymptotically fairly precise (at least when $R$ is even):

\proclaim{Theorem 2} If $\alpha$ has an expansion of period
$$ -R,R,-(R+1),(R+1),\tag 1.20 $$
then
$$ \rho(\alpha)=\cases {\displaystyle \frac{1}{4}\left( 1-\frac{3}{R}+\frac{4}{R^2}+O(R^{-3})\right),}  & \hbox{ if $R$ is even, } \\
{\displaystyle \frac{1}{4}\left( 1-\frac{3}{R}+\frac{5}{R^2}+O(R^{-3})\right),}  & \hbox{ if $R$ is odd. } \endcases \tag 1.21 $$
\endproclaim

Notice that Theorem 1 gives the absolute lower bound $\rho(\alpha)\geq 1/27.211...$ for all $\alpha$  (this being the best bound producable from  our choice $\gamma_*$).
However if we use what we shall call an {\it adjusted nearest integer expansion} rather than the nearest integer expansion 
then we can improve this bound slightly.

\proclaim{Theorem 3}
If (1.4) gives an `adjusted nearest integer expansion' for $\alpha$, as defined 
in Section 2 below,
then the corresponding $\gamma_{*}$ defined in (1.18) satisfies
$$ M(\alpha,\gamma_*)\geq \frac{(\sqrt{10}-3)(7-\sqrt{13})}{(31-2\sqrt{10}-3\sqrt{13})} =\frac{1}{25.159...}. \tag 1.22 $$
\endproclaim

It will be clear from the proof that we can come arbitrarily close 
to this bound by taking expansions
of
period 
$$ (3,)^k-2,2,-3, (-4,-2,2,2,)^k \tag 1.23 $$
with  $k$ suitably large (although presumably $\gamma_*$ does not give $\rho(\alpha)$).

Bounds of the form, 
$$ \rho (\alpha)\geq C_0 \tag 1.24 $$
for all $\alpha$ in $\R$,
are certainly not new (see  Davenport [\DavOne]). Davenport [\DavTwo]  gives $C_0= \frac{1}{128}$,
Ennola 
 [\Enn] improving this to 
$$C_0=\frac{1}{(16+6\sqrt{6})}=\frac{1}{30.24...}\tag 1.25 $$
(again best possible for his construction).
The bound 1/32 is given in [\RockBk].

The smallest known value of $\rho (\alpha)$ (giving an upper bound
on the largest $C_0$ possible in such a bound) still appears to be
Pitman's example [\Pitman]  
$$  \rho \left( \frac{\sqrt{3122285}-1097}{1094}\right)=\frac{547}{4\sqrt{3122285}}=\frac{1}{12.921...}.\tag1.26$$ 
This $\alpha$ has
a nearest integer expansion of  period 
$$3,-3,-3,-3,-2,3,-3,-2,\tag1.27 $$
 and an adjusted nearest integer expansion of period 
$$-4,-2,2,-4,-3,-2,2,-3.\tag1.28 $$

Davenport's and Ennola's bounds in fact  arise from the closely related
study of the inhomogeneous minima of real indefinite
quadratic forms (the connection is explained in Section 5).

\head 2. Notation \endhead

We recall some notations from [\paperone].  For an expansion of the form (1.4) we write $\alpha_i$ for the tail of the 
continued fraction 
$$ 
\alpha_{i}  :=[0;a_{i+1},\ve_{i+2}a_{i+2},\ve_{i+3}a_{i+3},\ldots ],  \tag 2.1 $$ 
and $\bar{\alpha}_{i}$
for the reversed expansion,
$$ \bar{\alpha}_{i}  : = [0;a_{i},\ve_{i}a_{i-1},\ve_{i-1}a_{i-2},\ldots ,\ve_2a_1 ]. \tag 2.2 $$ 
We shall primarily assume that (1.4) is the {\it nearest integer continued fraction expansion},  generated using the algorithm
$$\align  (a_0,\ve_1,\alpha_0) &   :=(0,1,\alpha), \\
 (a_{i+1},\ve_{i+2},\alpha_{i+1}) &  := \cases \left(\left\lfloor \alpha_{i}^{-1}\right\rfloor, 1,\left|\left|\alpha_{i}^{-1}\right|\right|\right),  & \hbox{ if $||\alpha_{i}^{-1}||=\{\alpha_{i}^{-1}\},$} \\
 \left(\left\lceil \alpha_{i}^{-1}\right\rceil, -1,\left|\left|\alpha_{i}^{-1}\right|\right|\right), & \hbox{ if $||\alpha_{i}^{-1}||=1-\{\alpha_{i}^{-1}\}$}. \endcases \;\; \;\; \tag 2.3 \endalign  $$
For the nearest integer continued fraction the `backwards' expansion will
be a {\it singular continued fraction expansion} (see Perron [\Perron, \S40])
and  we plainly have 
$$\alpha_i<1/2,\hskip3ex \bar{\alpha}_{i}\leq [0;2,\overline{-3}]=\frac{\sqrt{5}-1}{2}.\tag 2.4 $$ 
Writing  $p_i/q_i$ for the convergents
$$ \frac{p_i}{q_{i}}:=[a_0; \ve_1a_1,\ve_2a_2,...,\ve_ia_i], \tag2.5 $$
we define
$$ \delta_{i}:=-(-\ve_1)(-\ve_2)\cdots (-\ve_{i+1}), \;\;\; D_i:=\delta_{i}(q_i\alpha -p_i), \tag2.6 $$
and  recall the relations
$$ D_i=\alpha_0\alpha_1\cdots \alpha_i,\;\;\;\; q_i=(\bar{\alpha}_1\cdots \bar{\alpha}_{i})^{-1}.\tag2.7 $$

When $R=2$ there is an inherent asymmetry between the `forwards'
and  `backwards' expansions $\alpha_{i}$ and  $\bar{\alpha}_{i}$
in
that we can have $\ve_{i}=-1$, $a_i=2$ but not $a_i=2$, $\ve_{i+1}=-1$.
Thus there is some merit in using an {\it adjusted nearest integer expansion}
satisfying
\roster
\item"(i)" $a_i\geq 2$ with $a_i=2 \Rightarrow \ve_{i}$ or $\ve_{i+1}=1$.
\item"(ii)" If $a_i=2$, $\ve_{i}=-1$ then $1/\bar{\alpha}_{i-1}-1/\alpha_{i}\geq 1$. 
\item"(iii)"  If $a_i=2$, $\ve_{i+1}=-1$ then $1/{\alpha}_{i}-1/\bar{\alpha}_{i-1}> 1$. 
\item"(iv)" The expansion contains no block  $\cdots,2,-2,\cdots$ or $\cdots, 2, (-3,)^k-2,\cdots$.
\endroster
For these expansions we certainly have
$$ \alpha_i,\bar{\alpha}_{i} \leq [0;2,\overline{-3}]=\frac{\sqrt{5}-1}{2}.\tag 2.8 $$
Such an expansion can be obtained from the nearest integer continued fraction 
(and vice versa) by use of the relation
$$ [....,\ve a,-2,b,...]=[...,\ve (a-1),2,-(b+1),...] \tag 2.9$$
to replace any combinations $a_i=2,\ve_{i}=-1$ not satisfying (ii) by $a_i=2$, $\ve_{i+1}=-1$ satisfying (iii), applying
$$ \align 
[...,-3,-3,-2,b,...] & = [...,-3,-2,2,-(b+1),...] = [...,-2,2,-3,-(b+1),...] \tag2.10\endalign $$
as necessary to deal with any additional $a_i=2,\ve_{i}=-1$ not satisfying (ii) introduced
by this process by the presence of a preceding $-3$.
Blocks disallowed in (iv) will plainly not be introduced by this process.


For a real $\gamma$ in $[0,1]$ and expansion (1.4) 
we form the alpha expansion 
$$\gamma=\sum_{i=1}^{\infty} c_iD_{i-1} \tag2.11 $$ 
by setting $\gamma_{0}:=\gamma$ and 
$$ c_{i+1}:=\left\lfloor \frac{\gamma_{i}}{\alpha_{i}}\right\rfloor,\;\;\;  \gamma_{i+1}:=\left\{ \frac{\gamma_{i}}{\alpha_{i}}\right\}. \tag2.12 $$ 
Writing
$$ b_i:=\cases a_i, & \hbox{ if $\ve_{i+1}=1$,} \\
  a_i-1, &   \hbox{ if $\ve_{i+1}=-1$, } \endcases \tag2.13 $$
this  gives the unique expansion of the form (2.11) with coefficients satisfying

\roster
\item"(i)" $0\leq c_{i}\leq b_i  $,
\item"(ii)" If  $(c_i,\ve_{i+1})=(a_i,+1)$ then $c_{i+1}=0 $
\item"(iii)" If  $(c_i,\ve_{i+1})=(a_i-1,-1)$ with $(c_j,\ve_{j+1})=(a_j-2,-1)$ for any $i<
j<t$  then $c_{t}\leq b_t-1$.
If in addition $(c_t,\ve_{t+1})=(a_t-1,+1)$ then $c_{t+1}=0$.
\item"(iv)" The tail of the sequence does not consist solely of blocks of 
the form
{\;\;\;\;\;\roster
\item"\;\;\;(a)"
$ (c_i,\ve_{i+1})=(a_{i},+1),\;\;(c_{i+1},\ve_{i+1})=(0,\pm 1), $ 
\item"\;\;\;(b)" $(c_i,\ve_{i+1})=(a_i-1,-1),\;\;\;
(c_j,\ve_{j+1})=(a_j-2,-1),\;\; \hbox{ for 
all $j>i$},$
\item"\;\;\;(c)" $(c_i,\ve_{i+1})=(a_i-1,-1),\;\;(c_t,\ve_{t+1})=(a_t-1,+1),\;\;\;(c_{t+1},\ve
_{t+2})=(0,\pm 1), $
with $(c_j,\ve_{j+1})=(a_j-2,-1)$ for any $i<j<t$.
\endroster}
\endroster

From the alpha-expansion we extract the sequence of integers 
  such that
$$ c_i=\frac{1}{2}(a_i-1+\ve_{i}+t_{i}), \tag 2.14 $$
and form  
$d_{n}^+=d_{n}^+(\alpha,\gamma)$, $d_{n}^-=d_{n}^-(\alpha,\gamma)$ 
$$\align  d_{n}^- &  :=\left(\sum_{i=1}^{n} t_i\delta_{i-1}q_{i-1}\right)/ \delta_{n-1}q_{n} \\
 & =t_{n}\bar{\alpha}_{n}+ (-\ve_{n})t_{n-1}\bar{\alpha}_{n}\bar{\alpha}_{n-1}+ (-\ve_{n}) (-\ve_{n-1})t_{n-2}\bar{\alpha}_{n}\bar{\alpha}_{n-1}\bar{\alpha}_{n-2}+\cdots \\
    d_n^+ & :=\left(\sum_{i=n+1}^{\infty} t_i D_{i-1}\right)/D_{n-1}=t_{n+1}\alpha_n+t_{n+2}\alpha_n\alpha_{n+1}+t_{n+3}\alpha_n\alpha_{n+1}\alpha_{n+2}+\cdots . \tag 2.15 \endalign $$
We define
$$ S(\gamma):=\{i \; : \;  |d_i^-|\leq 1-\bar{\alpha}_{i} \hbox{ and }
-(1-\ve_{i+1}\alpha_i) \leq d_i^+ \leq (1+\ve_{i+1}\alpha_i) \}. \tag 2.16 $$

In Theorem 1 of  [\paperone] we showed:
\proclaim{Lemma 1} 
If $i$ is in $S(\gamma)$  for all but 
finitely many $i$ then
$$ M(\alpha,\gamma)=\liminf_{i\rightarrow \infty} \min \{ s_1(i),s_2(i),s_3(i),s_4(i)\}, $$
where
$$\align
s_1(i) & :=\frac{1}{4}(1-\bar{\alpha}_{i}+d_i^-)(1+\ve_{i+1}\alpha_i+d_i^+)/(1+\ve_{i+1}\bar{\alpha}_{i}\alpha_{i}), \\
s_3(i) & :=\frac{1}{4}(1-\bar{\alpha}_{i}-d_i^-)(1+\ve_{i+1}\alpha_i-d_i^+)/(1+\ve_{i+1}\bar{\alpha}_{i}\alpha_{i}), \\
s_2(i) & :=\frac{1}{4}(1+\bar{\alpha}_{i}-d_i^-)(1-\ve_{i+1}\alpha_i+d_i^+)/(1+\ve_{i+1}\bar{\alpha}_{i}\alpha_{i}), \\
s_4(i) & :=\frac{1}{4}(1+\bar{\alpha}_{i}+d_i^-)(1-\ve_{i+1}\alpha_i-d_i^+)/(1+\ve_{i+1}\bar{\alpha}_{i}\alpha_{i}). \tag2.17 \endalign $$
\endproclaim

We also gave  the upper bound:
\proclaim{Lemma 2} 
If $i$ is in $S(\gamma)$  for all but 
finitely many $i$ then
$$ M(\alpha,\gamma) \leq \liminf_{i\rightarrow \infty} \frac{1}{4} \left(\frac{a_i-|t_i|}{a_i+\ve_{i+1}\alpha_i+\ve_i\bar{\alpha}_{i-1}}\right), \tag 2.18 $$
with
$$ M(\alpha,\gamma) \leq \limsup_{i\rightarrow \infty} \frac{1}{4} \left(\frac{1}{a_i+\ve_{i+1}\alpha_i+\ve_i\bar{\alpha}_{i-1}}\right) \tag2.19 $$
otherwise.
\endproclaim
There is presumably a connection between Lemma 1 and the {\it divided cell 
algorithm} approach of Barnes [\BarnesTwo].

\head 3. Refined bounds \endhead

For an expansion (1.4) with $R=\liminf a_i$ we define
$$ R_*:=\cases R, & \hbox{ if $R$ is odd, } \\
   R+1, & \hbox{ if $R$ is even,} \endcases \;\;\;\;
R_{**}:=\cases R, & \hbox{ if $R$ is even, } \\
   R+1, & \hbox{ if $R$ is odd,} \endcases \;\;\;\; \tag 3.1 $$
and
$$ \theta_1:=\frac{1}{2}(\sqrt{R_*^2+4}-R_*),\;\;\; \theta_2:=\cases \frac{1}{2}(R_{**}-\sqrt{R_{**}^2-4}), & \hbox{ if $R\geq 3$}, \\
\frac{1}{2}, & \hbox{ if $R=2$.} \endcases \tag 3.2 $$
The lower bound (1.19) will follow from the precise bound:

\proclaim{Theorem 4}
If (1.4) gives the nearest integer expansion for $\alpha$,
then the corresponding $\gamma_{*}$ defined in (1.18) satisfies
$$ M(\alpha,\gamma_*)\geq M_1:=\frac{1}{4} \frac{\left(1-\frac{2}{R_*-\theta_2}\right)\left(1-\frac{1}{\max\{4,R_{**}\}}\left(1+\frac{\theta_1^2}{1-\theta_1}\right)\right)}{1+\frac{1}{\max\{4,R_{**}\}}\left(\theta_1-\frac{1}{R_*-\theta_2}\right)}. \tag 3.3$$

This bound is best possible in the sense that periodic $\alpha$ with period
$$ \align  (R_*,)^k R_{**},-R_*,  (-R_{**},)^k & \;\;\hbox{ if $R\geq 3$,} \\ (3,)^k4,-3, -2,A, &\;\; \hbox{ if $R=2$,}  \tag 3.4\endalign  $$
will have $M(\alpha,\gamma_*)\rightarrow M_1$ as $k$ (and $A$) 
$\rightarrow \infty$.
\endproclaim

We should note that  $\gamma_*$ will not give $\rho(\alpha)$ 
for the examples (3.4).
If $\ve_{i}=1$ for almost
all $i$ (i.e. the expansion agrees with the regular expansion from some point
onwards) then the lower bound (1.19) readily improves to 
$$ M(\alpha,\gamma_*)\geq \frac{1}{4}\left(1-\frac{2}{R}+O(\frac{1}{R^2})\right). \tag 3.5 $$
In this case it proves easier to remove the sign alternation
and consider
$$\gamma_{**}:=\sum_{i=1}^{\infty} \frac{1}{2}(a_i+\lambda_{i}'),\;\;\; \lambda_{i}':=\cases 0, & \hbox{ if $a_i$ is even,} \\ 1, & \hbox{ if $a_i$ is odd.}\endcases \tag 3.6  $$
This gives the optimal $\gamma$ when  $\alpha$ has
period  $+R^*$ (as does $\gamma_*$) and in a number of period two cases
(see [\papertwo])

\proclaim{ Theorem 5}

If (1.4) gives the regular continued fraction expansion of $\alpha$ with 
$R\geq 2$,                                   
then with $R_*$ and $\theta_1$ as in (3.1)  and (3.2)  above
$$\align  M(\alpha,\gamma_{**}) & \geq \frac{1}{4} \left(1-\frac{2}{R_*}\right)(1-2\theta_1^2-\theta_1^3)   \tag3.7  \\
  & \geq \frac{1}{4}\left(1-\frac{2}{R_*}-\frac{2}{{R_*}^2}+\frac{3}{R_*^3}\right). \endalign $$

\endproclaim

In particular when the regular continued fraction expansion
does not have infinitely many ones we obtain
$$ \rho(\alpha) \geq M(\alpha,\gamma_{**})\geq \frac{1}{2\sqrt{13}+8}=\frac{1}{15.211...}.\tag3.8 $$
Bound (3.7)  is sharp, achieved by taking an $\alpha$ with
infinitely many blocks 
$A,(R^*,)^k $
with $A$ and $k\rightarrow \infty$,
although $\gamma_{**}$ will not give $\rho(\alpha)$ in such cases.
Asymptotically though (3.7) is a good lower bound for $\rho (\alpha)$:

\proclaim{Theorem 6}
Suppose that the expansion of  $\alpha$ has period 
$$ R_{**},R_*,R_*,R_*,R_*. \tag 3.9 $$
Then
$$ \rho (\alpha)=\frac{1}{4}\left( 1-\frac{2}{R_{*}}-\frac{2}{R_{*}^2}+O(\frac{1}{R_{*}^3})\right). \tag 3.10 $$
\endproclaim

Setting
$$ \theta :=
\cases  \frac{1}{2}(R-\sqrt{R^2-4}), & \hbox{ if $R\geq 3$,} \\
   \frac{1}{2}, & \hbox{ if $R=2$,} \endcases  \tag3.11 $$
Theorem 1 also gives us  a bound of the form
$$ M(\alpha,\gamma_*)\geq \frac{1}{4}\left( \frac{(R-3)}{(R-2\theta)} +\frac{4}{R^3}\right), \;\;\; R\geq 4. \tag 3.12 $$ 
This immediately gives us a simple bound on the extent to which the integer
coefficients $c_i$ in the alpha-expansion  of optimal (in the sense of achieving $\rho (\alpha)$) or near optimal $\gamma$ can deviate
from $\frac{1}{2}(a_i-1+\ve_{i})$ :

\proclaim{Theorem 7}

Suppose  that  $R\geq 4$ is finite and that $\gamma$ has
$$ M(\alpha,\gamma)> M(\alpha,\gamma_*) -R_*^{-3}. \tag3.13  $$
If  (1.4) gives the nearest integer expansion of $\alpha$ then
$$ |t_i| <\frac{3a_{i}}{R}  \tag3.14 $$
for all but finitely many $i$.
In particular $a_k=R\geq 4$ implies that
$$ t_{k}=\cases \pm 1, & \hbox{ if $a_k$ is odd,} \\
    0 \hbox{ or } \pm 2,
& \hbox{ if $a_k$ is even,} \endcases \tag 3.15 $$ 
for almost all $k$.

If (1.4) gives the regular continued fraction expansion of $\alpha$ 
then
$$ |t_{i}|<\frac{2a_{i}}{R_*-1} \tag 3.16 $$
for almost all $i$ with 
$$ t_{k}=\cases \pm 1, & \hbox{ if $a_k$ is odd,} \\
    0, 
& \hbox{ if $a_k$ is even,} \endcases \tag 3.17 $$ 
when $a_{k}=R\geq 4$.

\endproclaim

The restriction $R\geq 4$ is actually needed here; for example a period two
alpha with $\ve_{2k}a_{2k}=-3,$ $\ve_{2k+1}a_{2k+1}=-4$
for all but finitely many $k$ has
$\rho(\alpha)$ achieved for $\gamma$ with $c_{2k}=2$ and $c_{2k+1}=0$ for all but finitely many $k$.
All possibilities in (3.15) and (3.17) can occur  (see examples in [\paperfour] and [\papertwo]).

\head 4. Proofs \endhead 

\head Proof of Theorem 1 \endhead

We show how Theorem 1 can be deduced from the lower bound $M_1$ of Theorem 4 (proved below).
We assume $R\geq 4$ (evaluating $M_1$ gives the bound for $R=2$ or 3).
Observe that
$$ 1+\frac{1}{\max\{4,R_{**}\}}\left(\theta_1-\frac{1}{R_*-\theta_2}\right)=1-\frac{\theta_1 (\theta_1+\theta_2)}{R_{**}(R_*-\theta_2)} < 1,$$
$$ 1-\frac{2}{R_*-\theta_2}=\cases 1-\frac{2}{R}-\frac{2}{R^3}+\frac{2\theta_2(R(1-\theta_2)-1)}{R^3(R-\theta_2)}, & \hbox{ if $R$ is odd,}\\
   1-\frac{2}{R}+\frac{2}{R^2}-\frac{4}{R^3}
+\frac{2(1-\theta_2)(3+\theta_2^2)}{R^3(R+1-\theta_2)}, & \hbox{ if $R$ is even,} \endcases $$
and
$$ 1-\frac{1}{\max\{4,R_{**}\}}\left(1+\frac{\theta_1^2}{1-\theta_1}\right)=\cases 1-\frac{1}{R+1}-\frac{1}{R^3}+\frac{\theta_1^2(R+2-\theta_1^2)}{R^3(R+1)}, & \hbox{ if $R$ is odd,}\\
   1-\frac{1}{R}-\frac{1}{R^3}+\frac{\theta_1^2(R(1+\theta_1)+1)}{R^3(1-\theta_1)}, & \hbox{ if $R$ is even.} \endcases $$
Hence if $R$ is odd
$$ \align M_1 & \geq \frac{1}{4} \left( 1-\frac{2}{R}-\frac{2}{R^3}\right)\left( 1-\frac{1}{R+1}-\frac{1}{R^3}\right) \\
 & =\frac{1}{4}\left( 1-\frac{3}{R}+\frac{3}{R^2}-\frac{6}{R^3}+\frac{(7R^3+2R^2+2R+2)}{R^6(R+1)}\right), \endalign $$
and if $R$ is even
$$ \align M_1 & \geq\frac{1}{4} \left( 1-\frac{2}{R}+\frac{2}{R^2}-\frac{4}{R^3}\right)\left(1-\frac{1}{R}-\frac{1}{R^3}\right)\\
  & =\frac{1}{4}\left(1-\frac{3}{R}+\frac{4}{R^2}-\frac{7}{R^3} +\frac{2(3R^2-R+2)}{R^6}\right),\endalign  $$
and the simplified lower bounds of Theorem 1 are plain. $\blacksquare$

\head Proof of Theorem 2  \endhead

We suppose that $\alpha$ has expansion of period 
$ -R_{**},R_{**},-R_{*},R_{*} $ and that $\gamma$ achieves $\rho (\alpha)$.
From Theorem 7  (proved below) we can assume that $t_{k}=0,\pm 2$ if $a_k=R_{**}$ and $t_{k}=\pm 1,\pm 3$ if $a_k=R_*$. Moreover we can eliminate $t_k=\pm 3$ since if $t_{k}=\pm 3$ infinitely often (2.17) gives
$$ M(\alpha,\gamma) \leq \frac{1}{4}\left(\frac{R_*-3}{R_*+O(R^{-2})}\right)=\frac{1}{4}\left(1-\frac{3}{R_*}+O(R^{-3})\right). $$
Now if $(t_k,t_{k+1})=(2,2)$ then
$$ \align s_{3}(k) & =\frac{1}{4}\left( 1-\frac{3}{R_{**}}+O(R^{-2})\right)\left(1-\frac{1}{R_{**}}+O(R^{-2})\right)/(1+O(R^{-2}))\\
  & =\frac{1}{4}\left(1-\frac{4}{R}+O(R^{-2})\right). \endalign $$
If $(t_k,t_{k+1})=(2,-2)$ then $s_{2}(k)$ gives the same result. 
Observing that once we have eliminated a block of $t_i$ from consideration
the same is true for its negative (since changing the signs of the $t_i$
merely interchanges $s_1(i)$ with $s_3(i)$ and $s_2(i)$ with $s_4(i)$) 
we can assume that
$(a_{k},a_{k+1})=(R_{**},R_{**})$ implies that $\pm (t_{k},t_{k+1})=(0,0),(2,0)$
or $(0,2)$. Now if $(t_{k},t_{k+1})=(0,2)$ then
$$\align  s_{3}(k+1) & \leq \frac{1}{4}\left( 1-\frac{3}{R_{**}}+O(R^{-3})\right)\left(1+\frac{1}{R_{*}^2}+O(R^{-3})\right)/\left(1-\frac{1}{R_{*}^2}+O(R^{-3})\right) \\
  & = \cases \frac{1}{4}\left( 1-\frac{3}{R}+\frac{2}{R^2}+O(R^{-3})\right), & \hbox{ if $R$ is even,} \\
  \frac{1}{4}\left( 1-\frac{3}{R}+\frac{5}{R^2}+O(R^{-3})\right), & \hbox{ if $R$ is odd,} \endcases \endalign $$
with $(2,0)$ being dealt with similarly using $s_3(k+1)$. Hence we assume
that $(t_{k},t_{k+1})=(0,0)$ when $(a_{k},a_{k+1})=(R_{**},R_{**})$.
Now if we have a block $\vec{t}=(t_{k},...,t_{k+5})=(0,0,1,1,0,0)$ then
$$\align  s_{3}(k+1) & \leq \frac{1}{4}\frac{\left( 1-\frac{1}{R_{**}}+O(R^{-3})\right)\left(1-\frac{2}{R_*}-\frac{1}{R_{*}^2}+O(R^{-3})\right)}{\left(1-\frac{1}{R_{*}^2}+O(R^{-3})\right)} \\
  & = \cases \frac{1}{4}\left( 1-\frac{3}{R}+\frac{4}{R^2}+O(R^{-3})\right), & \hbox{ if $R$ is even,} \\
  \frac{1}{4}\left( 1-\frac{3}{R}+\frac{3}{R^2}+O(R^{-3})\right), & \hbox{ if $R$ is odd,} \endcases \endalign $$
with $s_3(k+3)$ giving the same bound if $\vec{t}=(0,0,-1,1,0,0)$.
Since $\pm \vec{t}$ must be of one of these forms we gain the upper bound 
$(1.21)$ for $\rho (\alpha)$. For $R$ even the lower bound (1.19) shows that $M(\alpha,\gamma_*)$ achieves this, while for $R$ odd it is readily checked that 
the period four expansion $\overline{0,-2,1,1}$ for the $t_i$ achieves the upper bound (1.21). $\blacksquare$

\head  Proofs of The Lower Bounds in Theorems 3, 4 \& 5  \endhead

We suppose that we are working with the nearest integer expansion
or an adjusted  nearest integer expansion as appropriate. We assume that 
$R=\liminf a_i\geq 2$ is finite with  $a_i\geq R$ for all $i$ (from Theorem 1 of [\papertwo] changing a finite number of partial quotients $a_i$ or coefficients 
$c_i$ has no effect on the value of $M(\alpha,\gamma)$).
We assume that $\gamma_*$ has the expansion
(1.18), and when the nearest integer expansion is also a regular expansion
that $\gamma_{**}$ has expansion (3.6). From Lemmas 3 and 13 below  these expansions satisfy
the requirements of Lemma 1 and hence for Theorems 3, 4 and 5 it
will be enough to check that the functions $s_1(i),...,s_4(i)$
defined in  (2.17) are greater than the required bounds.

With  $R_*, R_{**}, \theta,\theta_1$ and $\theta_2$ as in (3.1), (3.2) 
and (3.11) above 
we define the parameters
$$ \align R_1  : =R_*- \theta, &\hskip3ex  R_2  : = \max\{2,  R_{**}-\theta\},\\
R_3  := R_*-\theta_2,   &\hskip3ex  R_4  : = \max\left\{ 4, R_{**}\right\}- \theta,\;\;\; R_5  : = R_{**}-\theta,\\
  \tag4.1\endalign $$
so that for a nearest integer expansion
$$ \alpha_{i}\leq \cases 1/R_1, & \hbox{ if $a_{i+1}$ is odd,} \\
1/R_2, & \hbox{ if $a_{i+1}$ is even,} \endcases  \;\;\;\; 
\bar{\alpha}_{i}\leq \cases 1/R_1, & \hbox{ if $a_{i}$ is odd,}\\
1/R_5, & \hbox{ if $a_{i}$ is even,} \endcases \tag4.2 $$
with $R_4$ replacing $R_5$ when $a_i\geq 4$. 

We give the following bounds on $d_i^+$ and $d_i^-$ for $\gamma_*$:

\proclaim{Lemma 3} If (1.4) is a nearest integer or adjusted nearest
integer expansion then
expansion (1.18) gives the valid alpha-expansion for $\gamma_*$, with all $i$ in $S(\gamma_*)$ (see (2.16)).

If (1.4) is a  nearest integer expansion then for $\gamma_*$
the  $d_i^+$, $d_i^-$ of (2.15)
satisfy
$$ |d_i^+|\leq \min\left\{ \alpha_{i}, \frac{1}{R_1}\right\},\;\;\; |d_i^-|\leq \frac{1}{1-\theta_1}\min\{ \theta_1,\bar{\alpha}_{i}\}. $$
If (1.4) is an adjusted nearest integer expansion with $R=2$ then $\gamma_*$
has
$$  |d_i^+|\leq \min\left\{ \alpha_{i}, \frac{2}{5}\right\},\;\;\; |d_i^-|\leq \min\left\{ \frac{1}{2},\frac{3}{2}\bar{\alpha}_{i}\right\}. $$
\endproclaim

\noindent
{\bf Proof:}
 Clearly $c_i=\frac{1}{2}(a_i-1+\ve_i+\lambda_i)$ has  $0\leq c_i <a_i$ with
$c_i=a_i-1$ implying $\ve_{i}=1$. Moreover  $c_i=a_i-1$ or $(c_i,\ve_{i})=(a_i-2,-1)$ can only occur when  $a_i=2$ or $a_i=3$ and $\lambda_i=1$. The sign
alternation in the $\lambda_i$ thus prevents $(c_i,\ve_{i+1})=(a_i-1,-1)$ 
being followed by $(c_j,\ve_{j+1})=(a_j-2,-1)$ for both $j=i+1$ and $i+2$.
Hence (1.18) gives a valid expansion.

Next we bound $d_i^+$. If $a_{i+j}$ is even for all $j\geq 1$ then plainly $d_i^+=0$,
so assume that $a_{i+J}$ is odd for some $J\geq 1$ with $a_{i+j}$ even for
any $1\leq j<J$. Then the alternation in the sign of the  $\lambda_i$ gives
$$ |d_i^+|=\left(\prod_{1\leq j < J} \alpha_{i+j-1}\right)|d_{i+J-1}^+|\leq \prod_{1\leq i\leq J}\alpha_{i+j-1}\leq \min\{\alpha_i, \alpha_{i+J-1}\} $$
and the bound follows from (4.2).

We assume first that we are using the nearest integer expansion and show
by induction that $\gamma_*$ satisfies $|d_i^-|\leq \theta_1/(1-\theta_1)$, the bound
$$ |d_i^-|\leq \bar{\alpha}_i\left(1+|d_{i-1}^-|\right) \tag 4.3 $$
then immediately giving $|d_i^-|\leq \bar{\alpha}_{i}/(1-\theta_1)$.
For $i=1$ we plainly have $|d_1^-|=0$ or $1/a_1\leq 1/R_*$ as $a_1$ is even or odd. If $a_i$ is even we have $|d_i^-|=\bar{\alpha}_{i}|d_{i-1}^-|<|d_{i-1}^-|$
so we assume that $a_{i}$ is odd. If $a_i>R_*$ or $a_i=R_*$, $\ve_{i}=1$, $\bar{\alpha}_{i-1}\geq \theta_1$ we have $\bar{\alpha}_{i}\leq \theta_1$ and (4.3)
gives
$$ |d_i^-|\leq \theta_1 \left( 1 +\theta_1/(1-\theta_1)\right). $$
If $a_i=R_*$, $\ve_{i}=1$, and $\bar{\alpha}_{i-1}\leq \theta_1$ then 
$\bar{\alpha}_{i}=1/(R_*+\bar{\alpha}_{i-1})$ and (4.3) gives
$$ |d_i^-|\leq \bar{\alpha}_{i} \left( 1 +\bar{\alpha}_{i-1}/(1-\theta_1)\right) <\theta_1/(1-\theta_1). $$
If $a_i=R_*$, $\ve_{i}=-1$ and $a_{i-1}$ is odd then the alternation in the
sign of the  $\lambda_i$ gives 
$$ |d_i^-|\leq \bar{\alpha}_{i} <\theta_1/(1-\theta_1), $$
since plainly $\bar{\alpha}_{i-1}<1-\theta_1$. 
If $a_i=R_*$, $\ve_{i}=-1$ and $a_{i-1}$ is even then
$$ |d_i^-|\leq \bar{\alpha}_{i}\left(1+\bar{\alpha}_{i-1}|d_{i-2}^-|\right) \leq \left(1+\bar{\alpha}_{i-1}\theta_1/(1-\theta_1)\right)/(R_*-\bar{\alpha}_{i-1}) <\theta_1/(1-\theta_1), $$
since $a_{i-1}\geq 4$ certainly ensures that $\bar{\alpha}_{i-1}<\frac{1}{2}(1-\theta_1)$. 
For the adjusted nearest integer expansion with $R=2$ we
only bother to obtain the less precise  bound $|d_i^-|<1/2$. 
The proof is almost identical, noting that if  $a_i$ is odd with 
$a_i\geq 5$  or 
$a_i=3$ and  $\ve_{i}=1$ then (4.3) gives
$$ |d_i^-|\leq \frac{1}{3}\left(1+\frac{1}{2}\right), $$
while  if  $a_i=3$ and  $\ve_{i}=-1$, then  $\bar{\alpha}_{i-1}<1/2$ gives
$$ |d_i^-|\leq \left(1+\frac{1}{2}\bar{\alpha}_{i-1}\right)/(3-\bar{\alpha}_{i-1}) <\frac{1}{2}. $$

Finally $i$ is in $S(\gamma_*)$ for all $i$ since 
$|d_i^-|\leq \frac{1}{2}\leq 1-\bar{\alpha}_i$ unless $\ve_i a_i=-2$ when
$|d_i^-|=\bar{\alpha}_i |d_{i-1}^-|\leq \frac{1}{2}\bar{\alpha}_{i}<1-\bar{\alpha}_i$, while $|d_i^+|\leq \frac{2}{5}<1-\alpha_i$ unless $a_{i+1}=2$ 
when $|d_{i}^+|=\alpha_i|d_{i+1}^+|<\frac{2}{5}\alpha_i<1-\alpha_i$.
 $\; \;\blacksquare$

\head Lemmas for Theorem 3 \endhead

We suppose that (1.4) gives an adjusted nearest integer expansion for $\alpha$
as defined in Section 2. We assume that $R=2$ (when $R\geq 3$ the result follows from Theorem 1). As remarked above we need only show that the $s_1(i),...,s_4(i)$ of
(2.17) are all at least $1/25.159...$. We break this proof down into a number of lemmas.

\proclaim{Lemma 4}
Suppose that (1.4) is an adjusted nearest integer expansion, as defined
in Section 2, with $\gamma_*$ as in (1.18). Then when $\ve_{i+1}=-1$ and  $a_{i+1}\geq 3$ 
$$B_1:=\frac{(1-\alpha_i+|d_i^+|)}{1-\bar{\alpha}_{i}\alpha_i}\geq \frac{\sqrt{10}+2}{\sqrt{10}+2 +2(1-\bar{\alpha}_{i})}, $$
with this bound achieved for $\alpha_i=\theta_0:=\frac{1}{3}(4-\sqrt{10})$,
this having period four expansion $\overline{-4,-2,2,2}$.

\endproclaim

\noindent
{\bf Proof. } We suppose that for some $J\geq 0$ we have
$$ \ve_{i+4j+2}a_{i+4j+1},...,\ve_{i+4j+5}a_{i+4j+4}=-4,2,2,-2, $$
for any $0\leq j<J$. So
$$ B_1=\frac{A_J-B_J\alpha_{i+4J}+|d_{i+4J}^+|}{C_J-D_J\alpha_{i+4J}} $$
where, writing $\mu:=19+6\sqrt{10},$ $\beta:=19-6\sqrt{10},$
$$ A_J:=A\mu^J+B\beta^J,\;\;\; C_J=C\mu^J+D\beta^J, $$
$$ B_J:=A\left(2-\frac{1}{2}\sqrt{10}\right)\mu^J+B\left(2+\frac{1}{2}\sqrt{10}\right)\beta^J,$$ 
$$D_{J}:=C\left(2-\frac{1}{2}\sqrt{10}\right)\mu^J+D\left(2+\frac{1}{2}\sqrt{10}\right)\beta^J,\;$$
with
$$ A:=\frac{1}{2\sqrt{10}}(\sqrt{10}+2),\;\; B:=\frac{1}{2\sqrt{10}}(\sqrt{10}-2),\;\;$$ 
$$ C:=\frac{1}{2\sqrt{10}}(\sqrt{10}+2+2(1-\bar{\alpha}_{i})),\;\; D:=\frac{1}{2\sqrt{10}}(\sqrt{10}-2-2(1-\bar{\alpha}_{i}).$$
Hence $B_1\rightarrow A/C$, the stated bound, as $J\rightarrow \infty$ with 
$$ \frac{A_J-B_J\alpha_{i+4J}+\kappa \alpha_{i+4J}}{C_J-D_J\alpha_{i+4J}} \geq A/C $$
when
$$ w:=(1-\bar{\alpha}_{i})(\theta_0-\alpha_{i+4J})\beta^J+(1-\theta_0\bar{\alpha}_{i})\kappa \alpha_{i+4J} \geq 0. $$
We can therefore assume that $\alpha_{i+4J}>\theta_0$ (since otherwise
$\kappa=0$  gives $w\geq 0$). If $a_{i+4J+1}=3$ then
$|d_{i+4J}^+|=\alpha_{i+4J}(1-|d_{i+4J+1}^+|)$ so we can take $\kappa=\frac{3}{5}$
and $w >0$ since $(1-\theta_0\bar{\alpha}_{i})>(1-\bar{\alpha}_{i})$
and $\alpha_{i+4J}<2/5$ gives $(\alpha_{i+4J}-\theta_0)<\frac{3}{5}\alpha_{i+4J}$. So we can assume that $a_{i+4J+1}=4$, $\ve_{i+4J+2}=-1$ and
(since $\alpha_{i+4J}>\theta_0$ implies $a_{i+4J+2}+\ve_{i+4J+3}\alpha_{i+4J+2}<\frac{1}{3}(4+\sqrt{10})$) also $a_{i+4J+2}=2,$ $\ve_{i+4j+3}=1$.
Now since $(4-\bar{\alpha}_{i+4J})-(a_{i+4J+3}+\ve_{i+4J+4}\alpha_{i+4J+3})\geq 1$ we must have $a_{i+4J+3}=2$ or 3 (with $\ve_{i+4J+4}=-1$ if $a_{i+4J+3}=3$).
If $a_{i+4J+3}=3$ then $|d_{i+4J}^+|=\alpha_{i+4J}\alpha_{i+4J+1}\alpha_{i+4J+2}(1-|d_{i+4J+3}^+|)=\alpha_{i+4J}(1-|d_{i+4J+3}^+|)/(7-2\alpha_{i+4J+3})$ and
we can take $\kappa =3/35$, giving $w >0$ since $\alpha_{i+4J}=(7-2\alpha_{i+4J+3})/(25-7\alpha_{i+4J+3})<7/25$ plainly gives $(\alpha_{i+4J}-\theta_0)<\frac{3}{35}\alpha_{i+4J}$. So we assume that $a_{i+4J+3}=2$ and the condition
$\alpha_{i+4J}>\theta_0$ gives $\ve_{i+4J+4}\alpha_{i+4J+3}>\frac{1}{2}(\sqrt{10}-2)$ forcing $\ve_{i+4J+4}=1,$ $a_{i+4J+4}=2$ and $\ve_{i+4J+5}=-1$. 
Thus the smallest value of $B_1$ occurs by letting $J=\infty$. $\blacksquare$

\proclaim{Lemma 5}
Suppose (1.4) is an adjusted nearest integer expansion
with $R=2$. If $\ve_{i-2}a_{i-2}=-2$ and $\pm d_{i-3}^+\leq 0$,  then $\gamma_*$
satisfies
$$ T_1:=\frac{(3-\bar{\alpha}_{i-3}\mp d_{i-3}^-)}{\left( \frac{1}{3}(11-\sqrt{10})-\bar{\alpha}_{i-3}\right)}\geq \left( \frac{3-\theta_1-\theta_1/(1-\theta_1)}{\frac{1}{3}(11-\sqrt{10})-\theta_1}\right), $$ 
the value when $\bar{\alpha}_{i-3}=\theta_1:=[0; \overline{3}]=\frac{1}{2}(\sqrt{13}-3).$
\endproclaim

\noindent
{\bf Proof. }
We suppose that for some $J\geq 0$ we have $\ve_{i-j-3}a_{i-j-3}=3$ for
any $0\leq j<J$. Then
$$ T_1=\frac{A_J+A_{J-1}\bar{\alpha}_{i-J-3}\pm (-1)^{J+1}d_{i-J-3}^-}{B_J+B_{J-1}\bar{\alpha}_{i-J-3}}, $$
where
$$ A_J:=A\theta_1^{-J}+B(-\theta_1)^{J}+\frac{1}{3}, \;\;\; B_J:=C\theta_1^{-J}+D(-\theta_1)^{J}, $$
with
$$ A=\frac{4}{3\sqrt{13}}(\sqrt{13}+2),\;\;\;B=\frac{4}{3\sqrt{13}}(\sqrt{13}-2),\;\;\;$$
$$ C=\frac{1}{3\sqrt{13}}((11-\sqrt{10})\theta_1^{-1}-3),\;\;\;D=\frac{1}{3\sqrt{13}}((11-\sqrt{10})\theta_1+3).$$
Hence $T_1\rightarrow A/C$, the stated bound, as $J\rightarrow \infty$ with
$$ \frac{A_J+(A_{J-1}-\kappa)\bar{\alpha}_{i-J-3}}{B_J+B_{J-1}\bar{\alpha}_{i-J-3}}>\frac{A}{C}, $$
when
$$ w_1:=\left(\frac{1}{3}(11-\sqrt{10})-\theta_1\right)\left(1+\bar{\alpha}_{i-J-3}(1-3\kappa)\right)> $$ 
$$ \hskip 10ex w_2:=\frac{4}{3}(5-\sqrt{10})
\left(\bar{\alpha}_{i-J-3}-\theta_1\right)(-\theta_1)^{J}. $$ 
We suppose that $\ve_{i-J-3}a_{i-J-3}\neq 3$. 
For $\bar{\alpha}_{i-J-3}=0$ we have $w_1=2.309...>0.741...\geq |w_{2}|,$
so it will be enough to check that $w_1>|w_2|$ for the largest 
appropriate value of $\bar{\alpha}_{i-J-3}$. If $a_{i-J-3}\geq 5$ then $\kappa=\frac{3}{2}$ and $\bar{\alpha}_{i-J-3}=\frac{2}{9}$ give $w_1\geq 0.513,$ $|w_2|<0.198$. If $a_{i-J-3}\geq 2$ is even then $\kappa=\frac{1}{2}$, $\bar{\alpha}_{i-J-3}=\frac{2}{3}$ give  $w_1\geq 1.539,$ $|w_2|<0.892$.
So we can assume that $a_{i-J-3}=3$, $\ve_{i-J-3}=-1$. If $a_{i-J-4}$ is odd
then $\pm (-1)^{J+1}d_{i-J-3}^->-\bar{\alpha}_{i-J-3}$ and taking $\kappa=1$, $\bar{\alpha}_{i-J-3}=\frac{2}{5}$ gives $w_1>0.461$, $|w_2|<0.239$. 
If $a_{i-J-4}\geq 4$ is even then  $\pm (-1)^{J+1}d_{i-J-3}^->-\bar{\alpha}_{i-J-3}(1+\frac{1}{2}\bar{\alpha}_{i-J-4})$ with $\bar{\alpha}_{i-J-4}<2/7$. Hence
taking $\kappa=\frac{8}{7}$ and $\bar{\alpha}_{i-J-3}=\frac{7}{19}$ gives 
$w_1>0.243$, $|w_2|<0.161$. 
So we can assume that $a_{i-J-4}=2$ and hence (from condition (iv) 
of an adjusted nearest integer expansion) $J\geq 1$  and $(3+\alpha_{i-J-3})-(a_{i-J-5}+\ve_{i-J-5}\bar{\alpha}_{i-J-6})\geq 1$ implies $a_{i-J-5}=2$ and $\ve_{i-J-5}\bar{\alpha}_{i-J-6}<\alpha_{i-J-3}<2/5$. So  $\pm (-1)^{J+1}d_{i-J-3}^->-\bar{\alpha}_{i-J-3}(1+\frac{1}{2}\bar{\alpha}_{i-J-4}\bar{\alpha}_{i-J-5})$ with $\bar{\alpha}_{i-J-4}\bar{\alpha}_{i-J-5}=(5+2\ve_{i-J-5}\bar{\alpha}_{i-J-6})^{-1}<\frac{1}{4}$ and 
$\bar{\alpha}_{i-J-3}=(5+2\ve_{i-J-5}\bar{\alpha}_{i-J-6})/(13+5\ve_{i-J-5}\bar{\alpha}_{i-J-6})<29/75$ and taking $\kappa=\frac{9}{8}$, $\bar{\alpha}_{i-J-3}=29/75$ gives $w_1>0.188$, $|w_2|<0.063$. Hence in all cases
we obtain $w>0$ and hence a larger contribution than taking $J=\infty$.
$\blacksquare$

\proclaim{Lemma 6}
Suppose that  (1.4) is an adjusted nearest integer expansion.
If $\ve_{i+1}a_{i+1}\neq 2$  then $\gamma_*$ satisfies
$$ s_2(i),s_4(i)=\frac{1}{4} (1+\bar{\alpha}_{i}\mp d_i^-)\frac{(1-\ve_{i+1}\alpha_i\pm d_i^+)}{(1+\ve_{i+1}\bar{\alpha}_{i}\alpha_i)} > \frac{1}{25.147...}. $$
\endproclaim

\noindent
{\bf Proof. } 
If $\ve_{i+1}=-1$ then $B_2:=(1-\ve_{i+1}\alpha_i\pm d_i^+)/(1+\ve_{i+1}\bar{\alpha}_{i}\alpha_i)>1$ and 
$$ s_2(i),s_4(i)= \frac{1}{4}(1+\bar{\alpha}_{i}\mp d_i^-)B_2 \geq \frac{1}{4}\left(1-\frac{1}{2}\bar{\alpha}_{i}\right)>\frac{1}{6}, $$
so we assume that $\ve_{i+1}=1$ and $a_{i+1}\geq 3$.  Hence
$$ B_2 \geq \frac{1-2\alpha_i}{1+\bar{\alpha}_{i}\alpha_{i}}>\frac{1}{5}\frac{1}{\left(1+\frac{2}{5}\bar{\alpha}_{i}\right)}, $$
and if $a_i$ is even we have
$$ s_2(i),s_4(i)\geq \frac{1}{4}\left(1+\frac{1}{2}\bar{\alpha}_{i}\right)\frac{1}{5}\frac{1}{(1+\frac{2}{5}\bar{\alpha}_i)}>\frac{1}{20}, $$
and if $a_i\geq 5$ is odd 
$$ s_2(i),s_4(i)\geq \frac{1}{4}\left(1-\frac{1}{2}\bar{\alpha}_{i}\right)\frac{1}{5}\frac{1}{(1+\frac{2}{5}\bar{\alpha}_i)}>\frac{1}{24.5}, $$
and we may assume that $a_i=3$. If $\alpha_i\leq \frac{1}{3}$ then 
$$ s_2(i),s_4(i)\geq \frac{1}{4}\left(1-\frac{1}{2}\bar{\alpha}_{i}\right)\frac{1}{3}\frac{1}{(1+\frac{1}{3}\bar{\alpha}_i)}>\frac{1}{17}, $$
so we may also assume that $a_{i+1}=3$, $\ve_{i+2}=-1$. Now if $\ve_{i}=-1$
and $a_{i-1}=2$ then 
$3+(3-\alpha_{i+1})^{-1}-\bar{\alpha}_{i-2}^{-1}\geq 1$, giving $\bar{\alpha}_{i}\leq [0;3,-2,2,(3-\alpha_{i+1})]$ with $\alpha_{i+1}<1/2$,  and in all cases $\bar{\alpha}_{i}<29/75.$ Likewise we obtain $\alpha_i<29/75$ and in this remaining case
$$ s_2(i),s_4(i)\geq \frac{1}{4}\left(1-\frac{1}{2}\bar{\alpha}_{i}\right)\frac{\frac{17}{75}}{(1+\frac{29}{75}\bar{\alpha}_i)}>\frac{1}{25.147...}. $$

\proclaim{Lemma 7}
Suppose that (1.4) is an adjusted nearest integer expansion. If $\pm d_i^+<0$
and $\bar{\alpha}_{i}<1/2$ and $\ve_{i+1}a_{i+1}\neq -3$ then $\gamma_*$ satisfies
$$ s_1(i),s_3(i)=\frac{1}{4}(1-\bar{\alpha}_{i}\pm d_i^-)\frac{(1+\ve_{i+1}\alpha_i\pm d_i^+)}{(1+\ve_{i+1}\bar{\alpha}_{i}\alpha_i)} >\frac{1}{24}. $$
\endproclaim

\noindent
{\bf Proof. }
Since $\pm d_i^+<0$  we have
$$ (1-\bar{\alpha}_{i}\pm d_i^-)\geq \cases 1-\frac{1}{2}\bar{\alpha}_i, & \hbox{ if $a_i$ is odd,} \\   1-\frac{3}{2}\bar{\alpha}_{i}, & \hbox{ if $a_i$ is even.}\endcases $$ 
If $\ve_{i+1}=1$ then 
$$ B_1:=\frac{(1+\ve_{i+1}\alpha_i\pm d_i^+)}{(1+\ve_{i+1}\bar{\alpha}_{i}\alpha_i)} > \frac{1}{1+\bar{\alpha}_{i}\alpha_i}>\frac{9}{13} $$
and
$$ s_1(i),s_3(i)=\frac{1}{4}(1-\bar{\alpha}_{i}\pm d_i^-)B_1 \geq \frac{9}{52}\left(1-\frac{3}{2}\bar{\alpha}_i\right)\geq \frac{1}{23.111...},$$
so we assume that $\ve_{i+1}=-1$. 
Hence 
$$ B_1\geq \cases \frac{1-7/5\alpha_{i}}{1-\bar{\alpha}_{i}\alpha_i}, & \hbox{ if $a_{i+1}$ is even,} \\ \frac{1-2\alpha_{i}}{1-\bar{\alpha}_{i}\alpha_i}, & \hbox{ if $a_{i+1}$ is odd,} \endcases $$
and if $a_i$ is odd then $\bar{\alpha}_i<2/5$ gives
$$ s_1(i),s_3(i) \geq \frac{1}{4}\left(1-\frac{1}{2}\bar{\alpha}_{i}\right)\frac{\frac{1}{5}}{1-\frac{2}{5}\bar{\alpha}_{i}} \geq \frac{1}{21}. $$
So we assume that $a_i$ is even. If $a_{i+1}\geq 4$ is even then 
$$ s_1(i),s_3(i) \geq \frac{1}{4}\left(1-\frac{3}{2}\bar{\alpha}_{i}\right)\frac{\frac{3}{5}}{1-\frac{2}{7}\bar{\alpha}_{i}} \geq \frac{1}{22.857...}, $$
and if $a_{i+1}=2$ then $a_i\geq 4$ gives $\bar{\alpha}_{i}<2/7$ and  
$$ s_1(i),s_3(i) \geq \frac{1}{4}\left(1-\frac{3}{2}\bar{\alpha}_{i}\right)\frac{\frac{3}{10}}{1-\frac{1}{2}\bar{\alpha}_{i}} \geq \frac{1}{20}, $$
so we assume that $a_{i+1}$ is odd.  If $\alpha_{i}<1/5$ then
$$ s_1(i),s_3(i) \geq \frac{1}{4}\left(1-\frac{3}{2}\bar{\alpha}_{i}\right)\frac{\frac{3}{5}}{1-\frac{1}{5}\bar{\alpha}_{i}} \geq \frac{1}{24}, $$
so we can assume that  $a_{i+1}=5$, $\ve_{i+2}=-1$. 
If $a_i=2$ then we must have $\bar{\alpha}_{i-1}>(4-\alpha_{i+1})^{-1}$ so that we can assume that $\bar{\alpha}_{i}<4/9$ and 
$$ s_1(i),s_3(i) \geq \frac{1}{4}\left(1-\frac{3}{2}\bar{\alpha}_{i}\right)\frac{\frac{5}{9}}{1-\frac{2}{9}\bar{\alpha}_{i}} \geq \frac{1}{19.466...}. $$

\proclaim{Lemma 8}
Suppose that (1.4) is an adjusted nearest integer expansion
with $R=2$. If $\pm d_i^+>0$
and $\bar{\alpha}_{i}<1/2$  then $\gamma_*$ satisfies
$$ \align s_1(i),s_3(i) & =\frac{1}{4}(1-\bar{\alpha}_{i}\pm d_i^-)\frac{(1+\ve_{i+1}\alpha_i\pm d_i^+)}{(1+\ve_{i+1}\bar{\alpha}_{i}\alpha_i)}\\
  &  \geq \frac{1}{4}(\sqrt{10}-3)\left( \frac{3-\theta_1-\theta_1/(1-\theta_1)}{\frac{1}{3}(11-\sqrt{10})-\theta_1}\right)=\frac{1}{25.159...},\endalign  $$
with asymptotic equality when
$...,\ve_{i}a_{i},|\ve_{i+1}a_{i+1},...=...,(3,)^k-2,2,-3,|(-4,-2,2,2,)^k...$
as $k\rightarrow \infty$. \endproclaim

\noindent
{\bf Proof. }  {\bf Step 1:} We first deal with $a_i\neq 3$.

  When $\pm d_i^+>0$ we have
$$ \align
  B_1:=\frac{(1+\ve_{i+1}\alpha_i\pm d_i^+)}{(1+\ve_{i+1}\bar{\alpha}_{i}\alpha_i)} & \geq \cases 1, & \hbox{ if $\ve_{i+1}=1$,} \\
\frac{1-\alpha_{i}}{1-\bar{\alpha}_{i}\alpha_{i}}, & \hbox{ if $\ve_{i+1}=-1$, } \endcases \\
  & \geq \frac{1}{2-\bar{\alpha}_{i}}. \endalign $$
Hence if $\bar{\alpha}_{i}\leq 2/7$ we have 
$$ s_1(i),s_3(i)=\frac{1}{4}(1-\bar{\alpha}_{i}\pm d_i^-)B_1\geq \frac{1}{4}\frac{\left(1-\frac{5}{2}\bar{\alpha}_{i}\right)}{\left(2-\bar{\alpha}_{i}\right)}\geq \frac{1}{24}. $$
If $a_i=2$ then 
$\bar{\alpha}_{i}<1/2$ gives
$$ s_1(i),s_3(i)\geq \frac{1}{4}\frac{(1-\frac{3}{2}\bar{\alpha}_{i})}{(2-\bar{\alpha}_{i})}\geq \frac{1}{24}, $$
so we assume that $a_i=3$ with $\bar{\alpha}_{i-1}<1/2$ if $\ve_{i}=1$.
Thus if $a_{i+1}=2,$ $\ve_{i+1}=-1$ we have $(3+\ve_{i}\bar{\alpha}_{i-1})-\alpha_{i+1}^{-1}\geq 1$ giving $\alpha_{i}\leq 2/5$ if $\ve_{i}=-1$ and $\alpha_i<5/12$ if $\ve_{i}=1$ and
$$ B_1 \geq \cases \frac{3}{5-2\bar{\alpha}_{i}}, & \hbox{ if $\ve_{i}=-1$,} \\
\frac{7}{12-5\bar{\alpha}_{i}}, & \hbox{ if $\ve_{i}=+1$.} \endcases $$

\noindent
{\bf Step 2:} Dealing with $\ve_ia_i=3$.

If $\ve_{i}=1$ then 
$$ s_1(i),s_3(i) \geq \frac{7}{4} \frac{(1-2\bar{\alpha}_{i}\mp d_{i-1}^-\bar{\alpha}_{i})}{(12-5\bar{\alpha}_{i})}=\frac{7}{4}\frac{(1+\bar{\alpha}_{i-1}\mp d_{i-1}^-)}{(31+12\bar{\alpha}_{i-1})}. $$
Hence if $a_{i-1}$ is even the bound $\pm d_{i-1}^-\geq -\frac{1}{2}\bar{\alpha}_{i-1}$ gives
$$  s_1(i),s_3(i) \geq \frac{7}{4} \frac{(1+\frac{1}{2}\bar{\alpha}_{i-1})}{(31+12\bar{\alpha}_{i-1})} > \frac{1}{17.714...} $$
and if $\bar{\alpha}_{i-1}<1/3$ then 
$$ s_1(i),s_3(i) \geq \frac{7}{4} \frac{(1-\frac{1}{2}\bar{\alpha}_{i-1})}{(31+12\bar{\alpha}_{i-1})} > \frac{1}{24}, $$
leaving only $a_{i-1}=3$, $\ve_{i-1}=-1$ in which case the observation that 
$\mp d_{i-1}^-\geq -\bar{\alpha}_{i-1}$ if $a_{i-2}$ is odd with $\mp d_{i-1}^-\geq -\bar{\alpha}_{i-1}-\frac{1}{2}\bar{\alpha}_{i-1}\bar{\alpha}_{i-2}\geq -\frac{5}{4}\bar{\alpha}_{i-1}$ if $a_{i-2}$ is even gives
$$ s_1(i),s_3(i) \geq \frac{7}{4} \frac{(1-\frac{1}{4}\bar{\alpha}_{i-1})}{(31+12\bar{\alpha}_{i-1})} > \frac{1}{22.730...}. $$
 
\noindent
{\bf Step 3:} Dealing with $\ve_ia_i=-3$ and $\ve_{i-1}a_{i-1}\neq 2$.

When  $\ve_{i}=-1$ we have 
$$ s_1(i),s_3(i) \geq \frac{3}{4} \frac{(1-2\bar{\alpha}_{i}\pm d_{i-1}^-\bar{\alpha}_{i})}{(5-2\bar{\alpha}_{i})}=\frac{3}{4}\frac{(1-\bar{\alpha}_{i-1}\pm d_{i-1}^-)}{(13-5\bar{\alpha}_{i-1})}. $$
If $a_{i-1}$ is odd then $\pm d_{i-1}^->\frac{1}{2}\bar{\alpha}_{i-1}$ and 
$\bar{\alpha}_{i-1}<2/5$ gives 
$$ s_1(i),s_3(i) \geq \frac{3}{4} \frac{(1-\frac{1}{2}\bar{\alpha}_{i-1})}{(13-5\bar{\alpha}_{i-1})} > \frac{1}{18.333...}. $$
If $a_{i-1}\geq 4$ is even then the bound 
$\pm d_{i-1}^-\geq -\frac{1}{2}\bar{\alpha}_{i-1}$ is enough to give
$$ s_1(i),s_3(i) \geq \frac{3}{4} \frac{(1-\frac{3}{2}\bar{\alpha}_{i-1})}{(13-5\bar{\alpha}_{i-1})} > \frac{1}{25.066...}, $$
when  $\bar{\alpha}_{i-1}\leq \frac{1}{4}$. When
$\ve_{i-1}a_{i-1}=-4$, the observation that $\pm d_{i-1}^->0$ if $a_{i-2}$ is odd and 
$\pm d_{i-1}^->-\frac{1}{2}\bar{\alpha}_{i-1}\bar{\alpha}_{i-2}$ 
if $a_{i-2}$ is even, with $\bar{\alpha}_{i-1}\leq \frac{2}{7}$ gives
$$ s_1(i),s_3(i) \geq \frac{3}{4} \frac{(1-\frac{5}{4}\bar{\alpha}_{i-1})}{(13-5\bar{\alpha}_{i-1})} > \frac{1}{24}. $$
So we are left only to deal with $\ve_ia_{i}=-3$, $\ve_{i-1}a_{i-1}=2$. 

\noindent
{\bf Step 4:} Dealing with $\ve_{i-1}a_{i-1}= 2$, $\ve_ia_i=-3$ and $\ve_{i+1}=1$.

Now if $\ve_{i+1}=1$ we obtain
$$ s_1(i),s_3(i)\geq \frac{1}{4}(1-2\bar{\alpha}_{i}\pm \bar{\alpha}_{i}d_{i-1}^-)\frac{(1+\alpha_i)}{(1+\bar{\alpha}_{i}\alpha_{i})}=\frac{1}{4}\frac{(1-\bar{\alpha}_{i-1}\pm d_{i-1}^-)(1+\alpha_{i})}{3-\bar{\alpha}_{i-1}+\alpha_{i}}. $$
Moreover since $(3+\alpha_i)-(a_{i-2}+\ve_{i-2}\bar{\alpha}_{i-3})> 1$
we must have $\bar{\alpha}_{i-1}<(2+\alpha_i)/(5+2\alpha_i)<8/19, $
and $a_{i-2}=2$ or $3$ with $a_{i-2}=2$ if $\alpha_{i}<\frac{1}{2}$.
Thus if $\alpha_i<1/2$ we have $\pm d_{i-1}^-\geq -\frac{1}{2}\bar{\alpha}_{i-1}\bar{\alpha}_{i-2}> -\frac{1}{3}\bar{\alpha}_{i-1}$ and
$$ s_1(i),s_3(i)\geq \frac{1}{4} \frac{(1-\frac{4}{3}\bar{\alpha}_{i-1})}{(3-\bar{\alpha}_{i-1})}\geq \frac{1}{23.52}, $$
and if $\alpha_i>1/2$
$$ s_1(i),s_3(i)\geq \frac{3}{8} \frac{(1-\frac{3}{2}\bar{\alpha}_{i-1})}{(\frac{7}{2}-\bar{\alpha}_{i-1})}\geq \frac{1}{22.285...}. $$

\noindent
{\bf Step 5:} Dealing with $\ve_{i-1}a_{i-1}= 2$, $\ve_ia_i=-3$ and $\ve_{i+1}=-1$.

Now if $(\ve_{i-1}a_{i-1},\ve_{i}a_i,\ve_{i+1}a_{i+1})=(2,-3,-a_{i+1})$ then $\ve_{i-2}a_{i-2}=-2$ (since $a_{i-2}+\ve_{i-2}\bar{\alpha}_{i-3}<2-\alpha_i$). From condition (iv) we must have $a_{i+1}\geq 3$. 
So we can assume that 
$ (\ve_{i-2}a_{i-2},\ve_{i-1}a_{i-1},\ve_{i}a_{i},\ve_{i+1}a_{i+1})=(-2,2,-3,-a_{i+1})$ with $a_{i+1}\geq 3$
and from Lemmas 4 and 5
$$ \align  s_1(i),s_3(i) & \geq \frac{1}{4}(\sqrt{10}+2)\frac{(1-2\bar{\alpha_i}\pm \bar{\alpha}_{i}d_{i-1}^-)}{(\sqrt{10}+4-2\bar{\alpha}_{i})}=\frac{1}{4}(\sqrt{10}-3)T_1\\
  & \geq \frac{1}{4}(\sqrt{10}-3)\left( \frac{3-\theta_1-\theta_1/(1-\theta_1)}{\frac{1}{3}(11-\sqrt{10})-\theta_1}\right), \endalign $$
the required lower bound, with asymptotic equality as claimed. $\blacksquare$

\head Proof of the Theorem 3 \endhead

From Lemma 6 we have $s_2(i),s_4(i)>1/25.148$ except when $\ve_{i+1}a_{i+1}=2$ in which case $s_2(i),s_4(i)=s_1(i+1),s_3(i+1)$. 

Hence it is enough to check the $s_1(i),s_3(i)$. 
If $\ve_ia_i=-2$ 
then $s_1(i),s_3(i)=s_1(i-1),s_3(i-1)$ 
with $a_{i-1}\neq 2$ so we
can assume that $\bar{\alpha}_{i}\leq \frac{1}{2}$.
Now if $\ve_{i+1}a_{i+1}=-3$ with $\pm d_i^+<0$ then 
$s_1(i),s_3(i)=s_1(i+1),s_3(i+1)$ with  
$\pm d_{i+1}^+>0$ and $\bar{\alpha}_{i+1}<\frac{1}{2}$. 
Lemmas 7 and 8 immediately give the required lower bound in the remaining cases. $\blacksquare$


\head Lemmas for Theorem 4 \endhead

For Theorem 4 we assume that (1.4) gives the nearest integer expansion
for $\alpha$.
We again break the proof down into manageable lemmas:

\proclaim{Lemma 9} Suppose that (1.4) gives the nearest integer expansion of $\alpha$. If $E_0\geq F_0\geq \theta_2$ then $\gamma_*$
satisfies 
$$ \frac{F_0-\alpha_{i+1}+|d_{i+1}^+|}{E_0-\alpha_{i+1}}\geq \frac{F_0-\theta_2}{E_0-\theta_2} $$
with equality iff $\ve_{i+j+2}a_{i+j+1}=-R_{**}$ for all $j\geq 1$.

\endproclaim
\noindent
{\bf Proof.}
Suppose that for some $J\geq 0$ we have $(a_{i+j+1},\ve_{i+j+2})=(R_{**},-1)$
for any $1\leq j\leq J$. Then 
$$ Q:=\frac{F_0-\alpha_{i+1}+|d_{i+1}^+|}{E_0-\alpha_{i+1}} =
\frac{F_J-F_{J-1}\alpha_{i+1+J}+|d_{i+1+J}^+|}{E_J-E_{J-1}\alpha_{i+1+J}} $$
with
$$ F_j:=A\theta_2^{-j}+B\theta_2^j, \;\;\;  E_j:=C\theta_2^{-j}+D\theta_2^j,$$
where, writing  $\Delta:=\sqrt{R_{**}^2-4}, $
$$ A=\frac{1}{\Delta}(F_0/\theta_2-1), \;\;\; B=\frac{1}{\Delta}(1-F_0\theta_2),$$
$$ C=\frac{1}{\Delta}(E_0/\theta_2-1),\;\;\; D=\frac{1}{\Delta}(1-E_0\theta_2).$$
Now as $J\rightarrow \infty$ we have $Q\rightarrow A/C$. 
Hence we assume that $(a_{i+J+2},\ve_{i+J+3})\neq (R_{**},-1)$. Observe that 
$$ \frac{F_J-(F_{J-1}-\kappa)\alpha}{E_{J}-E_{J-1}\alpha} >\frac{A}{C} $$
if and only if 
$$ w:=(E_0-F_0)\theta_2^J(\theta_2-\alpha)+(E_0-\theta_2)\kappa \alpha >0.$$
Now if $a_{J+i+2}>R_{**}$ or $a_{i+J+2}=R_{**}$ and $\ve_{i+J+3}=1$ then $\alpha_{i+J+1}<\theta_2$ and $\kappa=0$ gives $w >0$ in all three cases. If $a_{J+i+2}=R_*$ then
$|d_{i+J+1}^+|\geq \frac{1}{2}\alpha_{i+J+1}$ gives $\kappa =\frac{1}{2}$
and $w>0$ since $E_0-F_0<E_0-\theta_2$ and $(\alpha_{i+J+1}-\theta_2) <\frac{1}{2}\alpha_{i+J+1}$. 
Hence $Q\geq A/C$ with equality iff $J= \infty$. $\blacksquare$

\proclaim{Lemma 10}
With $R_3$, $\theta_2$ as in (4.1) and (3.2) define
$$ S_1:=\frac{1-2/R_3}{1-\bar{\alpha}_{i}/R_3},\;\;\;
 S_2:=\frac{1-\theta_2}{1-\theta_2\bar{\alpha}_{i}},\;\;\;
S_3:=\frac{1-2/R_3}{1+\bar{\alpha}_{i}/R_3},\;\;\;$$
and for $\gamma_*$ set
$$ B_1:= \frac{1+\ve_{i+1}\alpha_i\pm d_i^+}{1+\ve_{i+1}\bar{\alpha}_{i}\alpha_{i}}, \;\;\;  B_2:= \frac{1-\ve_{i+1}\alpha_i\pm d_i^+}{1+\ve_{i+1}\bar{\alpha}_{i}\alpha_{i}}. \;\;\; $$  
If (1.4) is the nearest integer expansion then $B_1$ and $B_2$ satisfy 
$$ B_1\geq \cases S_1, & \hbox{ if $\pm d_i^+ \leq 0,$} \\   
 S_2, & \hbox{ if $\pm d_i^+ > 0$,} \endcases \;\;\; B_2\geq S_3.  $$ 
Moreover $B_i\rightarrow S_1$ if 
$$\ve_{i+1}a_{i+1},\ve_{i+2}a_{i+2},\ve_{i+3}a_{i+3},...=\cases -R_*,(-R_{**},)^k...,  & \hbox{ if $R\geq 3$, } \\ -3,-2,A,..., & \hbox{ if $R=2$, } \endcases $$
as $k$ or $A \rightarrow \infty$. \endproclaim
 
\noindent
{\bf Proof.}
Suppose first that $\pm d_{i}^+\leq 0$. 
If $\ve_{i+1}=1$  then
$$ B_1\geq \frac{1}{1+\bar{\alpha}_{i}\alpha_{i}}\geq \frac{1}{1+\bar{\alpha}_{i}} >S_1. $$
If $\ve_{i+1}=-1$ and $a_{i+1}$ is even then 
$$ B_1\geq \frac{1-\left(1+1/R_1\right)/R_2}{1-\bar{\alpha}_{i}/R_2}>S_1. $$
If $\ve_{i+1}=-1$ and $a_{i+1}$ is odd then
$$ B_1=\frac{ 1-2\alpha_{i}+\alpha_{i}|d_{i+1}^+|}{1-\bar{\alpha}_{i}\alpha_i} $$
and when $R$ is even or  $\ve_{i+2}=1$ or $a_{i+1}>R_*$ the bound $B_1\geq S_1$
follows from $\alpha_i\leq 1/R_3$. So we can assume that $a_{i+1}=R_*=R$, $\ve_{i+2}=-1$ and by Lemma 9 
$$ B_1=\frac{(R_*-2)-\alpha_{i+1}+|d_{i+1}^+|}{(R_*-\bar{\alpha}_{i})-\alpha_{i+1}} \geq S_1, $$
with $B_1\rightarrow S_1$ in the cases claimed.

Similarly when $\pm d_i^->0$ we have $B_1\geq (1+\alpha_i)/(1+\bar{\alpha}_{i}\alpha_i)>1$ if $\ve_{i+1}=1$ with   
$$B_1\geq (1-\alpha_i)/(1-\bar{\alpha}_{i}\alpha_i)>S_2$$ 
if $\ve_{i+1}=-1$ and  $R$ is even or $a_{i+1}>R_{**}$ or $a_{i+1}=R_{**}$ and 
$\ve_{i+2}=1$. For $\ve_{i+1}=-1$, $R=R_*$ and $a_{i+1}$ odd we have
$$ B_1 \geq \frac{1-\alpha_i/R_1}{1-\bar{\alpha}_{i}\alpha_i} >\frac{1-1/R_1^2}{1-\bar{\alpha}_{i}/R_1} >S_2. $$
Thus we may assume when $\pm d_i^+ > 0$  that $a_{i+1}=R_{**}=R+1$, $\ve_{i+2}=-1$ and by Lemma 9
$$ B_1=\frac{(R_{**}-1)-\alpha_{i+1}+|d_{i+1}^+|}{(R_{**}-\bar{\alpha}_{i})-\alpha_{i+1}}\geq S_2.  $$
Similarly for $B_2$ we have $B_2\geq 1$ if $\ve_{i+1}=-1$. If $\ve_{i+1}=1$ then
$$ B_2 \geq \frac{1-\alpha_{i}(1+1/R_1)}{1+\bar{\alpha}_{i}\alpha_{i}} \geq \frac{1-(1+1/R_1)/R_2}{1+\bar{\alpha}_{i}/R_2} >S_3 $$
when $a_{i+1}$ is even, with
$$ B_2\geq \frac{1-2\alpha_{i}+\alpha_{i}|d_{i+1}^+|}{1+\bar{\alpha}_{i}\alpha_{i}} $$
when $a_{i+1}$ is odd, with $\alpha_i<1/R_3$ giving $B_2\geq S_3$  unless $R$ is odd, $a_{i+1}=R_*$, $\ve_{i+2}=-1$ in which case by Lemma 9
$$ B_2=\frac{(R_*-2)-\alpha_{i+1}+|d_{i+1}^+|}{(R_*+\bar{\alpha}_{i})-\alpha_{i+1}} \geq S_3. \; \blacksquare $$

\proclaim{Lemma 11} Suppose that (1.4) is the nearest integer expansion
and define $d:=\max\{ 4,R_{**}\}$. 
For $(A_0,B_0,\delta_0)$ equal to
\roster
\item"(i)" $(d-1,d-1/R_3,1)$,
\item"(ii)" $(R_*-2,R_*-\theta_2,1)$,
\item"(iii)"$(R_*,R_*+1/R_3,-1)$,
\endroster
define  
$$ H:= \frac{ (A_0+\theta_1)-(1+\theta_1)/R_*}{B_0+\theta_1 }. $$
Then for $\gamma_*$ when $\pm \delta_0d_{i-1}^+\leq 0$  we have
$$ \frac{A_0+\ve_{i}\bar{\alpha}_{i-1}\mp \delta_0\ve_{i}d_{i-1}^-}{B_0+\ve_{i}\bar{\alpha}_{i-1}}\geq \cases \min\{ H, 4M_1/(1-\theta_2)\}, & \hbox{ in case (ii) when $R\leq 4$,} \\ H, & \hbox{ otherwise,} \endcases $$
with the expression tending to $H$ if $\ve_{j+1}a_j=+R_*$ for all $j<i$ as $i\rightarrow \infty$.
\endproclaim

\noindent
{\bf Proof. } Assume that for some $J\geq 0$ we 
have $(a_{i-j},\ve_{i-j+1})=(R_*,1)$ for any $1\leq j\leq J$. Then we can write
$$  T  :=\frac{A_0+\ve_{i}\bar{\alpha}_{i-1}\mp \delta_0\ve_{i}d_{i-1}^-}{B_0+\ve_{i}\bar{\alpha}_{i-1}}
   = \frac{ A_{J}+A_{J-1}\ve_{i-J}\bar{\alpha}_{i-J-1}\pm (-1)^{J+1}\delta_0\ve_{i-J}d_{i-J-1}^-}{B_{J}+B_{J-1}\ve_{i-J}\bar{\alpha}_{i-J-1}}, $$
where 
$$ A_j:=A\alpha^j+B\beta^j+1/R_*, \;\;\; B_j:=C\alpha^j+D\beta^j, $$
with
$ \alpha:=\theta_1^{-1}$, $\beta:=-\theta_1,$
and, writing $\Delta':=\sqrt{R_*^2+4},$
$$ A:=((A_0-\beta)-(1-\beta)/R_*)\frac{\alpha}{\Delta'}, \;\;\; B:=-((A_0-\alpha)-(1-\alpha)/R_*)\frac{\beta}{\Delta'}, $$
$$ C:=(B_0-\beta)\frac{\alpha}{\Delta'}, \;\;\; D:=-(B_0-\alpha)\frac{\beta}{\Delta'}. $$
As $J\rightarrow \infty$ we have $T\rightarrow A/C=H$. So we assume that $(a_{i-J-1},\ve_{i-J})\neq (R_*,1)$. 
We observe that 
$$ T'=\frac{A_J+(A_{J-1}-\kappa)\ve \bar{\alpha}}{B_J+B_{J-1}\ve\bar{\alpha}} >\frac{A}{C} $$
if and only if
$$ w := (R_*(A_0-B_0)+(B_0-1))\beta^J(\theta_1 -\ve \bar{\alpha}) + (\theta_1 +B_0)(1+\ve\bar{\alpha}(1-R_*\kappa))>0. \tag 4.4 $$ 
Now if $\ve_{i-J}=1$ and $a_{i-J-1}$ is even then $\pm (-1)^{J+1}\delta_0d_{i-J-1}^-\geq -\bar{\alpha}_{i-J-1}\theta_1/(1-\theta_1)$ and we obtain $T>T'$ with
$$ \ve =1,\;\;\; \kappa =\theta_1/(1-\theta_1),\;\;\; \bar{\alpha}_{i-J-1}\leq 1/(R_{**}-\theta). \tag 4.5 $$
If $\ve_{i-J}=1$ and $a_{i-J-1}>R_*$ is odd then 
$\pm (-1)^{J+1}\delta_0d_{i-J-1}^-\geq -(1+\theta_1/(1-\theta_1))\bar{\alpha}_{i-J-1}$ and we obtain $T'$ with
$$ \ve =1,\;\;\; \kappa =1/(1-\theta_1),\;\;\; \bar{\alpha}_{i-J-1}\leq \left(R_{*}+\frac{3}{2}\right)^{-1}. \tag 4.6 $$
If $\ve_{i-J}=-1$ and $a_{i-J-1}$ is odd then $\pm (-1)^{J+1}\delta_0\ve_{i-J}
d_{i-J-1}^-\geq (1-\theta_1/(1-\theta_1))\bar{\alpha}_{i-J-1}$ and we obtain $T'$ with
$$ \ve =-1,\;\;\; \kappa =1-\theta_1/(1-\theta_1),\;\;\; \bar{\alpha}_{i-J-1}\leq \left(R_{*}-\frac{1}{2}\right)^{-1}. \tag 4.7 $$
If $\ve_{i-J}=-1$ and $a_{i-J-1}\geq 4$ is even then $\pm (-1)^{J+1}\delta_0
\ve_{i-J}d_{i-J-1}^-\geq -\theta_1/(1-\theta_1)\bar{\alpha}_{i-J-1}$ and we obtain $T'$ with
$$ \ve =-1,\;\;\; \kappa =-\theta_1/(1-\theta_1),\;\;\; \bar{\alpha}_{i-J-1}\leq 
\left(a_{i-J-1}-\frac{1}{2}\right)^{-1}. \tag 4.8 $$
Now since 
$$\frac{|R_*(A_0-B_0)+B_0-1|}{(\theta_1+B_0)}\leq \frac{(R_*+1)-(R_*-1)\theta_2}{(R_*+\theta_1-\theta_2)} $$ 
for the three choices of $(A_0,B_0)$ 
we will certainly have $w >0$ in all cases where
$$ w_1 :=|\theta_1-\ve \bar{\alpha}_{i-J-1}|\theta_1^J((R_*+1)-(R_*-1)\theta_2) \tag 4.9 $$
$$\hskip 12ex \leq w_2:=(R_*+\theta_1-\theta_2)(1+\ve\bar{\alpha}_{i-J-1}(1-R_*\kappa)). $$ 
Now substituting  $\bar{\alpha}_{i-J-1}=0$ we  have 
$$ w_1\leq \theta_1 (R_*+1)<\frac{1}{3}(R_*+1)< R_*-\frac{1}{2}\leq w_2$$
so for (4.5)--(4.8) it will be enough to check that (4.9) holds for
the largest value of $\bar{\alpha}_{i-J-1}$. For (4.5) we have
$$ w_1< \frac{1}{2} (R_*+1-2\theta_2) \leq \frac{2}{3}(R_*-\theta_2)\leq (R_*+\theta_1-\theta_2)(1-\bar{\alpha}_{i-J-1}\theta_1)=w_2. $$
For (4.6) we have
$$ w_1 = \theta_1^{J+1}\bar{\alpha}(3/2-\theta_1) (R_*+1-(R_*-1)\theta_2)\leq \frac{1}{2}\bar{\alpha}\left(3/2-\theta_1\right) (R_*+1-2\theta_2), $$
while
$w_2 =\bar{\alpha} (3/2-\theta_1)(R_*+\theta_1-\theta_2).$
For (4.7) 
$$ w_1 \leq \frac{4}{5}(R_*+1-2\theta_2)\leq \frac{16}{15}(R_*-\theta_2), $$
while
$$ w_2=\left(2-\frac{(3/2+\theta_1)}{\left(R_*-1/2\right)}\right)(R_*+\theta_1-\theta_2)\geq \frac{19}{15} (R_*-\theta_2). $$
For (4.8) when $a_{i-J-1}\geq 6$ we have
$$ w_1\leq \frac{2}{11}\left( 1+\frac{11}{2}\theta_1\right)(R_*+1-(R_*-1)\theta_2), $$
$$ w_2=(1-\bar{\alpha}_{i-J-1}(2+\theta_1))(R_*+\theta_1-\theta_2)\geq \frac{2}{11}\left(\frac{7}{2}-\theta_1\right) (R_*+\theta_1-\theta_2), $$
and $w>0$ (using 
$w_1 \leq \frac{2}{5}(R_*+1-4\theta_2)$ and  $w_2\geq \frac{3}{5}(R_*-\theta_2)$ when $R_*\geq 5$ and checking numerically when $R_*=3$).
If $J\geq 1$ and $a_{i-J-1}=4$  then $w_1\leq \frac{1}{5}(R_*+1-2\theta_2)$ and $w_2\geq \frac{2}{7}(R_*-\theta_2)$ gives the required result. One readily checks numerically that $w>0$ still holds for $J=0$, $a_{i-1}=4$, $R=2,3,4$ 
in case (i) and (iii). Hence it remains only to
deal with case (ii) when $a_{i-1}=4$, $\ve_{i}=-1$ (the minimum occuring in the other cases by letting $J\rightarrow \infty$). 
When  $R=2,3$ and $a_{i-1}=4$, $\ve_{i}=-1$ we obtain
$$ T =\frac{1}{3-\theta_2} \left(\frac{3+\ve_{i-1}\bar{\alpha}_{i-2}\mp \ve_{i-1}d_{i-2}^-}{(4-1/R_3)+\ve_{i-1}\bar{\alpha}_{i-2}}\right) \geq
 \frac{1}{3-\theta_2}\left( \frac{3+\theta_1-(1+\theta_1)/3}{4-1/R_3+\theta_1}\right)=\frac{4M_1}{1-\theta_2}$$
using case (i). 
For $R=4$ and $\ve_{i}=-1$, $a_{i-1}=4$ we have
$$ T=\frac{1}{5-\theta_2} \left(\frac{11+3 \ve_{i-1}\bar{\alpha}_{i-2}\mp \ve_{i-1}d_{i-2}^-}{(4-1/R_3)+\ve_{i-1}\bar{\alpha}_{i-2}}\right) $$
and the estimates $\mp\ve_{i-1} d_{i-2}^-\geq -\frac{5}{4}\bar{\alpha}_{i-2}$ if $\ve_{i-1}=1$ and $\mp \ve_{i-1} d_{i-2}^-\geq-\frac{1}{4}\bar{\alpha}_{i-2}$ if $\ve_{i-1}=-1$
(since $\pm d_{i-2}^-\geq \frac{3}{4}\bar{\alpha}_{i-2}$ if $a_{i-2}$ is odd)
with $\bar{\alpha}_{i-2}<2/7$ are enough to give 
$$ \align \frac{1}{4}(1-\theta_2)T & \geq \frac{1}{4}\frac{(1-\theta_2)}{(5-\theta_2)}\min \left\{ \frac{11+\frac{7}{4}\bar{\alpha}_{i-2}}{(4-1/R_3)+\bar{\alpha}_{i-2}},
\frac{11-\frac{13}{4}\bar{\alpha}_{i-2}}{(4-1/R_3)-\bar{\alpha}_{i-2}}\right\} \\
  & \geq \frac{1}{9.160...}>M_1=\frac{1}{9.337...}.\;\;\;\blacksquare
\endalign  $$

\proclaim{Lemma 12}
With $\theta_1$, $\theta_2$, $M_1$ and $R_3$ as in (3.2), (3.3) and (4.1), define
$ d:=\max\{ 4, R_{**}\}, $
and
$$ U_1:=\frac{1-\frac{1}{(d+\theta_1)}\left(1+\frac{\theta_1}{1-\theta_1}\right)}{1-\frac{1}{R_{3}(d+\theta_1)}},\;\; U_2:=\frac{1-\theta_1-\frac{\theta_1}{1-\theta_1}}{1-\theta_1\theta_2},\;\; U_3:=\frac{1-\theta_1^2/(1-\theta_1)}{1+\theta_1/R_3}.  $$
For $\gamma_*$ set
$$ T_1:=\frac{1-\bar{\alpha}_{i}\pm d_{i}^-}{1-\bar{\alpha}_{i}/R_3},\;\;T_2:=\frac{1-\bar{\alpha}_{i}\pm d_i^-}{1-\bar{\alpha}_{i}\theta_2},\;\;T_3:=\frac{1+\bar{\alpha}_{i}\mp d_i^-}{1+\bar{\alpha}_{i}/R_3}. $$
Suppose that  (1.4) is the nearest integer expansion. Then
$ T_3\geq U_3. $

\noindent
When $\pm d_i^+\leq 0$ and $a_i\geq 3$
we have $  T_1\geq U_1, $
with $T_1\rightarrow U_1$ for $\ve_ia_{i},\ve_{i-1}a_{i-1},\ve_{i-2}a_{i-2},... =d,(R_*,)^k...$ as $k\rightarrow \infty$.

\noindent
When $\pm d_i^+\geq 0$ and $a_i\geq 3$ we have
$$ T_2\geq \cases U_2, & \hbox{ if $R\geq 5$, } \\
  \min\{ U_2, 4M_1/(1-\theta_2)\}, & \hbox{ if $R\leq 4$.} \endcases $$
\endproclaim

\noindent
{\bf Proof. }
Suppose first that $\pm d_i^+<0$ and $a_i\geq 3$. If $a_i$ is odd then 
$\pm d_i^-\geq (1-\theta_1/(1-\theta_1))\bar{\alpha}_{i}$ and
$$ T_1\geq \frac{1-\frac{1}{2}\bar{\alpha}_{i}}{1-\bar{\alpha}_{i}/R_3}\geq \frac{1-\frac{1}{2}1/R_1}{1-1/(R_3R_1)} >U_1, $$
and if $a_i>d+2$ is even then $\bar{\alpha}_{i}\leq 1/(d+2-\theta)$ 
gives 
$$ T_1 \geq \frac{1-(1+\theta_1/(1-\theta_1))\bar{\alpha}_{i}}{1-\bar{\alpha}_{i}/R_3} >U_1.$$  
So we assume that $a_i=d$ and hence from Lemma 11
$$ T_1=\frac{(d-1)+\ve_{i}\bar{\alpha}_{i-1}\mp \ve_{i}d_{i-1}^-}{(d-1/R_3)+\ve_{i}\bar{\alpha}_{i-1}} \geq U_1.  $$
Similarly if $\pm d_i^+>0$ and $a_i\geq 4$ is even then
$$ T_2 \geq \frac{1-(1+\theta_1/(1-\theta_1))\bar{\alpha}_{i}}{1-\bar{\alpha}_{i}\theta_2}\geq \frac{1-1/((1-\theta_1)R_4)}{1-\theta_2/R_4} >U_2, $$
and if $a_i>R_*$ is odd then $\bar{\alpha}_{i}<\theta_1$ and 
$$ T_2> \frac{1-(2+\theta_1/(1-\theta_1))\bar{\alpha}_{i}}{1-\bar{\alpha}_{i}\theta_2} >U_2. $$
So we may assume that $a_i=R_*$ and from Lemma 11
$$ T_2=\frac{(R_{*}-2)+\ve_{i}\bar{\alpha}_{i-1}\mp \ve_{i}d_{i-1}^-}{(R_*-\theta_2)+\ve_{i}\bar{\alpha}_{i-1}} \geq \min\{ U_2, 4M_1/(1-\theta_2).  $$
Likewise for $T_3$ observe that if $a_i$ is even or $\mp d_i^->0$ then 
$T_3>(1+\frac{1}{2}\bar{\alpha}_{i})/(1+\bar{\alpha}_{i}/R_3)>1$ and if $a_i>R_*$
is odd then $\bar{\alpha}_{i}<\theta_1$ and  $T_3>(1-\bar{\alpha}_{i}\theta_1/(1-\theta_1))/(1+\bar{\alpha}_{i}/R_3)>S_3$. So we assume that $a_i=R_*$ with $\mp d_i^-<0$ and by Lemma 11
$$ T_3=\frac{R_* +\ve_{i}\bar{\alpha}_{i-1}\pm \ve_{i}d_{i-1}^-}{(R_*+1/R_3)+\ve_{i}\bar{\alpha}_{i-1}}\geq U_3. \; \blacksquare $$

\head Proof of Theorem 4 \endhead

We assume that we are working with the nearest integer expansion 
and show that the $s_1(i),...,s_4(i)$ of (2.17)  are always greater than the
claimed bound $M_1$.

Now from Lemma 10 and Lemma 12
$$\align  \frac{1}{4} (1+\bar{\alpha}_{i}\mp d_{i}^-)\frac{(1-\ve_{i+1}\alpha_i\pm d_i^+)}{(1+\ve_{i+1}\bar{\alpha}_{i}\alpha_i)} & \geq \frac{1}{4}\left(1-\frac{2}{R_3}\right)\frac{ (1+\bar{\alpha}_{i}\mp d_{i}^-)}{(1+\bar{\alpha}_i/R_3)} \\ & \geq  \frac{1}{4}\left(1-\frac{2}{R_3}\right)
\frac{\left(1-\frac{\theta_1^2}{1-\theta_1}\right)}{\left(1+\frac{\theta_1}{R_{3}}\right)}. \endalign $$
Similarly if $\pm d_i^+\leq 0$ and $a_i\geq 3$, then
$$ \align  \frac{1}{4} (1-\bar{\alpha}_{i}\pm d_{i}^-)\frac{(1+\ve_{i+1}\alpha_i\pm d_i^+)}{(1+\ve_{i+1}\bar{\alpha}_{i}\alpha_i)} & \geq \frac{1}{4}\left(1-\frac{2}{R_3}\right)\frac{ (1-\bar{\alpha}_{i}\pm d_{i}^-)}{(1-\bar{\alpha}_i/R_3)} \\ 
 & \geq \frac{1}{4}\left(1-\frac{2}{R_{3}}\right)\frac{\left( 1-\frac{1}{(1-\theta_1)(\max\{4,R_{**}\}+\theta_1)}\right)}{\left(1-\frac{1}{R_3(\max\{4,R_{**}\}+\theta_1)}\right)}=M_1,  \endalign $$
with asymptotic equality when $...,\ve_ia_i,|\ve_{i+1}a_{i+1},...=...,(R_*,)^k,R_{**},|-R_*,(-R_{**},)^k,...$ when $R\geq 3$ and $...,(3,)^k 4,|-3,-2,A,...$ if $R=2$, 
as $k$ (and $A$) $\rightarrow \infty$.

If $\pm d_i^+> 0$ and $a_i\geq 3$, then
$$  \align  \frac{1}{4}(1-\bar{\alpha}_{i}\pm d_{i}^-)\frac{(1+\ve_{i+1}\alpha_i\pm d_i^+)}{(1+\ve_{i+1}\bar{\alpha}_{i}\alpha_i)} & \geq \frac{1}{4}\left(1-\theta_2\right)\frac{ (1-\bar{\alpha}_{i}\pm d_{i}^-)}{(1-\bar{\alpha}_i\theta_2)} \\ &  \geq \cases M_2, & \hbox{ if $R\geq 5$,} \\
   \min\{M_1,M_2\}, & \hbox{ if  $R\leq 4$.} \endcases  \endalign$$
with
$$ M_2:=\frac{1}{4}\frac{(1-\theta_2)}{(1-\theta_1\theta_2)}\left(1-\theta_1-\frac{1}{R_*}(1+\theta_1)\right). $$
The minimum of these bounds is $M_1$.

Hence it remains only to check $s_1(i),s_3(i)$ when $a_i=2$.
If $a_i=2$ and $\ve_{i}=-1$ then
$$ s_1(i)=s_1(i-1), \;\;\; s_3(i)=s_3(i-1), \;\;\; a_{i-1}\geq 3. $$
If $a_i=2$ and $\ve_{i+1}=\ve_{i}=1$ we have $\frac{(1+\ve_{i+1}\alpha_i\pm d_i^+)}{1+\ve_{i+1}\bar{\alpha}_{i}\alpha_i}\geq (1+\frac{2}{5}\bar{\alpha}_{i})^{-1}$ and $\pm d_{i-1}^-\geq -\frac{3}{2}\bar{\alpha}_{i-1}$, $\bar{\alpha}_{i-1}\leq (\sqrt{5}-1)/2$ giving
$$ s_1(i),s_3(i) \geq \frac{1}{4}\frac{(1-\bar{\alpha}_{i}\mp \bar{\alpha}_{i}d_{i-1}^-)}{\left(1+\frac{2}{5}\bar{\alpha}_{i}\right)} =\frac{1}{4}\frac{(1+\bar{\alpha}_{i-1}\mp d_{i-1}^-)}{\left(\frac{12}{5}+\bar{\alpha}_{i-1}\right)} \geq \frac{1}{4}\frac{(1-\frac{1}{2}\bar{\alpha}_{i-1})}{\left(\frac{12}{5}+\bar{\alpha}_{i-1}\right)}\geq \frac{1}{17.5}.\; \blacksquare $$

\head  Lemmas for  Theorem 5  \endhead
 
For Theorem 5 we assume that the nearest integer expansion has $\ve_i=1$ for
all $i$.

\proclaim{Lemma 13}
If (1.4) is a nearest integer expansion with all $\ve_i=1$ then (3.6) gives
a valid expansion for $\gamma_{**}$, with all $i$ in  $S(\gamma_{**})$ and
$$ 0<d_i^+\leq \frac{1}{R_*}+\frac{1}{R_*^2},\;\;\; |d_i^-|\leq \frac{1}{R_*}. $$
\endproclaim

\noindent
{\bf Proof:}
We suppose that $\ve_i=1$ for all $i$.
Plainly $c_i=\frac{1}{2}(a_i+\lambda_i')<a_i$ and expansion (3.6)
is a valid alpha-expansion. Suppose that  $a_{i+J}$ is odd for some $J\geq 1$, with $a_{i+j}$ even for any
$1\leq i<J$, then, since trivially $\lambda_i'\leq a_i/R_*$,
$$ 0\leq d_i^+D_{i-1} \leq \frac{1}{R_*}\sum_{j\geq i+J} a_jD_{j-1}=\frac{1}{R_*} D_{i+J-2}(1+\alpha_{i+J-1})\leq \frac{1}{R_*}\left(1+\frac{1}{R_*}\right)D_{i-1}, $$
with $d_i^+=0$ if all $a_j$, $j>i$, are even. Similarly 
$$ -\frac{1}{R_*}q_{i-1}\leq  -\frac{1}{R_*}\sum_{j\geq 0} a_{i-1-2j}q_{i-2j-2}\leq  q_id_i^-\leq \frac{1}{R_*} \sum_{j\geq 0} a_{i-2j}q_{i-2j-1}\leq \frac{1}{R_*}q_{i}. $$
Plainly $|d_{i}^{\pm}|<\frac{1}{2}\leq 1-\bar{\alpha}_{i},1-\alpha_i$ so $i$ is in $S(\gamma_{**})$ for all $i$.
 $\; \;\blacksquare$

\proclaim{Lemma 14}
Suppose that the nearest integer expansion (1.4)
is also a regular expansion. Then for $\gamma_{**}$
$$ L:=\frac{(1+\alpha_{i}-|d_i^+|)}{1+\alpha_{i}/R_{*}}\geq 1-2\theta_1^2-\theta_1^3, $$
with equality if $a_j=R_*$ for all $j>i$.

\endproclaim

\noindent
{\bf Proof. }
Suppose that for some $J\geq 0$ we have $a_{i+j}=R^*$ for any $0<j\leq J$
and $a_{i+J+1}\neq R_{*}$. Then
$$ L=\frac{A_{J}+A_{J-1}\alpha_{i+J} -|d_{i+J}^+|}{C_{J}+C_{J-1}\alpha_{i+J}} $$
where
$$ A_{J}=A\theta_1^{-J} -B(-\theta_1)^J+\frac{1}{R_*},\;\;\; C_{J}=C\theta_1^{-J} -D(-\theta_1)^J, $$
with
$$ A=\left(1-\frac{1}{R_*}\right)\frac{1+\theta_1^{-1}}{\sqrt{R_*^2+4}},\;\;\; 
 B=\left(1-\frac{1}{R_*}\right)\frac{1-\theta_1}{\sqrt{R_*^2+4}},\;\;\;$$ 
$$  C=\frac{1/R_*+1/\theta_1}{\sqrt{R_*^2+4}},\;\;\; 
 D=\frac{1/R_*-\theta_1}{\sqrt{R_*^2+4}},   $$
with
$ L\rightarrow A/C$, as $J\rightarrow \infty$.
Now 
$$ \frac{A_{J}+(A_{J-1}-\kappa)\alpha}{C_{J}+C_{J-1}\alpha} >A/ C $$
if and only if 
$$ w_1:=(1+\alpha (1-\kappa R_*)> w_2:=(R_*-1)^2(-\theta_1)^{J+1}(\alpha-\theta_1). $$
This is plainly true when $\alpha=0$ so we need only examine the largest
appropriate value of $\alpha_{i+J}$.
Now if $a_{i+J+1}\geq R_*+2$ is odd then 
$|d_{i+J}^+|=\alpha_{i+J}(1+|d_{i+J+1}^+|)$  with $\alpha_{i+J}<1/(R_*+2)$ 
and taking
$$ \kappa =1+\frac{1}{R_*}+\frac{1}{R_*^2}, \;\;\; \alpha =\frac{1}{R_*+2}, $$
we have
$$ |w_2|=\frac{\theta_1^{2+J}}{R_*+2}(2-\theta_1)(R_*-1)^2 <\frac{2\left(1-1/R_*\right)^2}{R_*+2}< \frac{2-1/R_*}{R_*+2} =w_1. $$
Similarly when $a_{i+J+1}$ is even we have $|d_{i+J}^+|=\alpha_{i+J}|d_{i+J+1}^+|$
and taking
$$ \kappa =\frac{1}{R_*}+\frac{1}{R_*^2},\;\;\;\; \alpha=\frac{1}{R_{**}}, $$
gives
$$ |w_2|< \theta_1^2(R_*-1)^2 <\left(1-\frac{1}{R_*}\right)^2
< \left(1-\frac{1}{R_{*}R_{**}}\right)=w_1. $$
Hence in all cases $L\geq A/C$. $\blacksquare$ 

\head Proof of Theorem 5 \endhead

We again need only show the desired lower bound for the $s_1(i),...,s_4(i)$
of (2.17).

We first deal with $s_2(i)$ and $s_4(i)$:

Suppose first that  $a_i$ is even, or $a_i$ odd and $\mp d_i^-\geq 0$.
Then  
$$ \frac{(1+\bar{\alpha}_{i}\mp d_i^-)}{1+\bar{\alpha}_{i}\alpha_{i}} \geq \frac{1+\bar{\alpha}_{i} \left(1-\frac{1}{R_*}\right)}{1+\bar{\alpha}_{i}\alpha_{i}}>1. $$
Now except when $a_{i+1}=3$ and $\pm d_{i}^+\leq 0$ we have
$$ 1-\alpha_{i}\pm d_i^+ \geq \cases 1-\frac{1}{R_{**}}\left(1+\frac{1}{R_*}+\frac{1}{R_*^2}\right), & \hbox{ if $a_{i+1}$ is even,} \\
 1-\frac{1}{\max\{R_*,5\}}-\left(\frac{1}{R_*}+\frac{1}{R_*^2}\right), & \hbox{ if $a_{i+1}$ is odd,} \endcases $$
giving
$$ s_2(i),s_4(i)=\frac{1}{4}\frac{(1+\bar{\alpha}_{i}\mp d_i^-)}{(1+\bar{\alpha}_{i}\alpha_{i})} (1-\alpha_{i}\pm d_i^+)  \geq \cases \frac{1}{4}\left(1-\frac{2}{R_{*}}-\frac{1}{R_*^2}\right),  &  \hbox{ if $R_*\geq 5$,} \\ 
\frac{1}{11.25}, & \hbox{ if $R_*=3$.} \endcases \tag 4.10 $$
If $a_{i+1}=3$ and $\pm d_i^+\leq 0$ then $s_2(i),s_4(i)=s_{1}(i+1),s_3(i+1)$
(to be dealt with below).

Suppose now that $a_i$ is odd and $\mp d_i^-\leq 0$.  
Then  $\pm d_i^+\geq 0$ with
$\pm d_i^+ =\alpha_i+\alpha_i|d_{i+1}^+|$ if $a_{i+1}$ is odd and
$$ \frac{1-\alpha_i\pm d_i^+}{1+\bar{\alpha}_{i}\alpha_i} \geq \cases  \frac{1-1/R_{**}}{1+\bar{\alpha}_{i}/R_{**}}, & \hbox{ if $a_{i+1}$ is even,} \\
\frac{1}{1+\bar{\alpha}_{i}/R_{*}}, & \hbox{ if $a_{i+1}$ is odd,} \endcases $$
while $\mp d_i^- =-\bar{\alpha}_{i}\pm \bar{\alpha}_{i}d_{i-1}^-$ gives
$$ \frac{1+\bar{\alpha}_{i}-|d_i^-|}{1+\bar{\alpha}_{i}/R_{**}} \geq \left(1-\frac{1}{R_*^2}\right)/\left(1+\frac{1}{R_*R_{**}}\right), $$
and
$$ s_2(i),s_4(i) \geq \frac{1}{4} \left(1-\frac{1}{R_{**}}\right)
\left(1-\frac{1}{R_*^2}\right)/\left(1+\frac{1}{R_*R_{**}}\right), $$
a bound larger than (4.10).

It remains to deal with $s_1(i)$, $s_3(i)$. 

If $a_i$ is odd then $|d_i^-|\leq 1/R_*$ gives 
$$ \frac{1-\bar{\alpha}_{i}\pm d_i^-}{1+\bar{\alpha_{i}}\alpha_{i}} > \frac{1-2/R_*}{1+\alpha_{i}/R_*}.$$
If $a_i$ is even then $|d_i^-|=\bar{\alpha}_{i}|d_{i-1}^-|\leq \bar{\alpha}_{i}/R_* $ gives
$$ \frac{1-\bar{\alpha}_{i}\pm d_i^-}{1+\bar{\alpha}_{i}\alpha_{i}}\geq \frac{1-\frac{1}{a_{i}}\left(1+\frac{1}{R_*}\right)}{1+\alpha_i/a_i}> \frac{1-2/R_*}{1+\alpha_{i}/R_*}, $$
as long as $a_i\geq 4$. If $a_i=2$ the bound $\frac{1}{3}(1+\frac{1}{3}\alpha_i)^{-1}$ still holds, the improved estimates $\bar{\alpha}_{i}\leq (2+\frac{1}{6})^{-1}$ if $a_{i-1}\leq 5$ and 
$|d_{i-i}^-|\leq \bar{\alpha}_{i-1}(1+|d_{i-2}^-|)\leq \frac{4}{3}\bar{\alpha}_{i-1}\leq \frac{2}{9}$ if $a_{i-1}\geq 6$ giving
$\frac{5}{13}(1+\frac{6}{13}\alpha_i)^{-1}$ and
$\frac{7}{18}(1+\frac{1}{2}\alpha_i)^{-1}$ respectively.
Hence, applying Lemma 14, 
$$ s_1(i),s_3(i) =\frac{1}{4}\frac{(1-\bar{\alpha}_{i}\pm d_i^-)}{(1+\bar{\alpha_{i}}\alpha_{i})} (1+\alpha_i\pm d_{i}^+) \geq \frac{1}{4} \left(1-\frac{2}{R_{*}}\right) L\geq
\frac{1}{4} \left(1-\frac{2}{R_{*}}\right) ( 1-2\theta_1^2-\theta_1^3). $$
It  is readily checked that this is smaller than (4.10),
giving the lower bound claimed.
 $\blacksquare$ 

\head Proof of Theorem 6 \endhead
Suppose that $\alpha$ has period  
$R_{**},R_*,R_*,R_*,R_*, $ and that $\gamma$ achieves $\rho (\alpha)$.
 From Theorem 7
 we can assume that 
$t_k=\pm 1$ if $a_k=R_*$ and $t_k=0,\pm 2$ if $a_k=R_{**}$. We first rule out
$t_{k}=\pm 2$ (this is immediate from (3.16) 
if $R$ is even as observed in (3.17)). Suppose that $a_{k}=R_{**}$
and $t_k=2$. If $t_{k+1}=1$ then
$$ s_{3}(k)= \frac{1}{4}\left(1-\frac{3}{R_{**}}+O(\frac{1}{R^2})\right)\frac{\left(1+O(R^{-2})\right)}{\left(1+O(R^{-2})\right)} =\frac{1}{4}\left(1-\frac{3}{R_*}+O(R^{-2})\right), $$
and if $t_{k+1}=-1$ then
$$ \align s_{2}(k) & =\frac{1}{4}\left(1-\frac{1}{R_{**}}+O(R^{-2})\right)\left(1-\frac{2}{R_{*}}+O(R^{-2})\right)/(1+O(R^{-2})) \\
  & =\frac{1}{4}\left(1-\frac{3}{R_*}+O(R^{-2})\right). \endalign $$
Similarly $t_k=-2$ is dealt with using $s_1(k)$ and $s_4(k)$. So we assume
that $t_k=0$ if $a_{k}=R_{**}$. Suppose that $a_{k}=R_{**}$ and set $\vec{t}=(t_{k},t_{k+1},...,t_{k+5})$. If $\vec{t}=(0,1,1,1,*,0)$ 
then 
$$\align  s_3(k+1) &  =\frac{1}{4}\frac{\left(1-\frac{2}{R_{*}}+O(R^{-3})\right)\left(1-\frac{1}{R_{*}^2}+O(R^{-3})\right)}{\left(1+\frac{1}{R_*^2}+O(R^{-3})\right)} \\
  & =\frac{1}{4}\left(1-\frac{2}{R_*}-\frac{2}{R_*^2}+O(R_*^{-3})\right).\endalign  $$
If $\vec{t}=(0,*,-1,1,-1,0)$ then we obtain the same expression for $s_2(k+3)$.
If $\vec{t}=(0,1,-1,-1,*,0)$ 
then 
$$\align  s_2(k+1) &  =\frac{1}{4}\frac{\left(1-\frac{2}{R_{*}}-\frac{1}{R_*^2}+O(R^{-3})\right)\left(1+O(R^{-3})\right)}{\left(1+\frac{1}{R_*^2}+O(R^{-3})\right)} \\
  & =\frac{1}{4}\left(1-\frac{2}{R_*}-\frac{2}{R_*^2}+O(R_*^{-3})\right).\endalign  $$ 
Similarly for $\vec{t}=(0,*,-1,1,1,0)$ using $s_3(k+3)$.
We can also immediately  dispose of $-\vec{t}$ of these forms
(by observing that changing the signs of the $t_k$ merely interchanging the roles of $s_1$ and $s_3$, and $s_2$ and $s_4$).
Since this includes all possible choices of $\vec{t}$ we can do no better
than (3.10), this  being achievable with $\gamma_{**}$. $\blacksquare$

\head Proof of Theorem 7 \endhead

Suppose that $R\geq 4$ and that $\gamma$ has 
$$M(\alpha,\gamma)> \frac{R-3}{R-2\theta}.$$
Noting that 
$$\frac{1}{a_j+\ve_j\alpha_j+\ve_{j-1}\bar{\alpha}_{j-1}} \leq \frac{1}{R-2\theta}, $$
we can assume from (2.19) that $i$ is in $S(\gamma)$ for almost all $i$.
Hence from (2.18)  we have 
$$ M(\alpha,\gamma) \leq \frac{1}{4} \liminf_{i\rightarrow \infty} \frac{(a_i-|t_{i}|)}{(a_i+\ve_i\alpha_{i}+\ve_{i-1}\bar{\alpha}_{i-1})} \leq \frac{1}{4}\frac{\liminf_{i\rightarrow \infty}(1-|t_i|/a_i)}{(1-2\theta/R)}, $$
and we must have $|t_i|< 3a_i/R$ for all but finitely many $i$. 
For the regular continued fraction we have 
 $M(\alpha,\gamma_{**})\geq \frac{1}{4}\left(1-2/(R_*-1)+4/R_*^3\right)$
for  $R\geq 4$, 
and the second bound follows similarly. 
$\blacksquare$

\head 5. Lower bound problems for quadratic forms \endhead

 For a real indefinite  form
$$ f(x,y)=ax^2+bxy+cy^2=a(x-\alpha y)(x-\bar{\alpha}y)\in \R[x,y], \tag 6.1 $$
with
$$ D(f):=b^2-4ac > 0, \tag 6.2 $$
and a point $P=(r,s)$ in $\R^2\setminus \Z^2$ one defines 
$$ \align M(f,P) & :=\inf_{(x,y)\in \Z^2} |f(x+r,y+s)|,\\
 M^*(f,P) & :=\liminf_{\Sb (x,y)\in \Z^2\\ \max\{ |x|,|y|\}\rightarrow \infty\endSb} |f(x+r,y+s)|. \tag 6.3 \endalign $$
The {\it inhomogeneous minima}
is then
$$  M(f)   :=\sup_{P} M(f,P). \tag 6.4 $$
It is straightforward to see that 
$$ M^*(f,P)=\sqrt{D(f)} \min\{ M (\alpha,\gamma_P),M(\bar{\alpha},\bar{\gamma}_{P})\},\;\;\;\gamma_P=r-s\alpha,\;\;\bar{\gamma}_P=r-s\bar{\alpha}.\tag 6.5$$
Hence, defining a related constant $M^*(f)\geq M(f)$ by
$$    M^*(f)   :=\sup_{P} M^*(f,P), \tag  6.6 $$
a bound of the form
$$ M^*(f)\geq C_0 \sqrt{D(f)} \tag 6.7 $$
for all forms $f$, is equivalent to (1.24). Bound (6.7)  certainly follows from the bound
$$  M(f)\geq C_0 \sqrt{D(f)}  \tag 6.8 $$
actually proved by Davenport and Ennola.
We should note that in the case of historical interest the $f$ are
irreducible forms in $\Z[x,y]$ and one has $M(f)/\sqrt{D(f)}=\rho(\alpha)$.
Bounds of the form (6.8) originally rose out of the study of norm-Euclidean
real quadratic fields (see Barnes [1] for details).

\enddocument
\end